\documentclass[11pt]{article}
\usepackage[utf8]{inputenc}
\usepackage[a4paper,top=2.5cm,bottom=2cm,left=1.5cm,right=1.5cm,marginparwidth=1.75cm]{geometry}

\usepackage[numbers]{natbib}
\usepackage{graphicx}
\usepackage{amsmath}
\usepackage{amssymb}
\usepackage{amsthm}
\usepackage[utf8]{inputenc}
\usepackage{dsfont}
\usepackage{textcomp}
\usepackage{url}
\usepackage{sidecap}
\usepackage{mathtools}
\usepackage{chngcntr}
\usepackage{MnSymbol}
\usepackage{wasysym}
\usepackage{graphicx}
\usepackage{setspace}
\usepackage{framed}
\usepackage{xcolor} 
\usepackage{transparent}
\usepackage{blindtext}
\usepackage{float}
\usepackage{subcaption}
\usepackage{authblk}
\usepackage{bbm}
\providecommand{\keywords}[1]{\textbf{\textbf{Keywords}} #1}

\newcommand{\R}{\mathbb{R}}

\newcommand{\1}{\mathbbm{1}}

\newcommand{\E}{\mathbb{E}}

\theoremstyle{remark}
\newtheorem*{remark}{Remark}

\title{\textbf{From interacting agents to \\   density-based modeling with stochastic PDEs}}

\author[a,b,*]{Luzie Helfmann}
\author[b]{Nata\v sa Djurdjevac Conrad}
\author[a,c]{Ana Djurdjevac}
\author[b]{Stefanie Winkelmann}
\author[a,b]{Christof Sch\" utte}

\affil[a]{Freie Universit\" at Berlin, Institut f\"ur Mathematik, Berlin, Germany}
\affil[b]{Zuse Institute Berlin, Berlin, Germany}
\affil[c]{Technische Universit\" at Berlin, Institut f\"ur Mathematik, Berlin, Germany}
\affil[*]{Corresponding author: luzie.helfmann@fu-berlin.de}

\date{November 2020}

\begin{document}
\newgeometry{
  left=1.9cm,
  right=1.9cm,
  top=2.5cm,
  bottom=3cm,
  bindingoffset=5mm
}

\maketitle

\begin{abstract}
Many real-world processes can naturally be modeled as systems of interacting agents. However, the long-term simulation of such agent-based models is often intractable when the system becomes too large.  In this paper, starting from a stochastic spatio-temporal agent-based model (ABM), we present a reduced model in terms of stochastic PDEs that describes the evolution of agent number densities for large populations  while  retaining the inherent model stochasticity. We discuss the algorithmic details of both approaches; regarding the SPDE model, we apply Finite Element discretization in space which not only ensures efficient simulation but also serves as a regularization of the SPDE. Illustrative examples for the spreading of an innovation among agents are given and used for comparing  ABM and SPDE models.
\end{abstract}

 
\section{Introduction} 
Modeling real-world dynamics is of great importance for the understanding, prediction and manipulation of the dynamical process of interest. Mathematically, one is therefore interested in the model analysis, inference of its parameters, as well as the accurate simulation and control of the dynamics. 
The modeling of realistic processes is nowadays often based on formulations in terms of discrete and autonomous entities (e.g., humans, institutions, companies), so-called agents, that move in a given environment and interact with each other according to a set of rules.
Agent-based models (ABMs) are very flexible in their description, often they just consist of a set of behavioural rules for agents encoded in computer programs~\cite{bonabeau2002agent,macy2002factors,galan2009errors}. Due to their versatility and the ease of incorporating data and assumptions into the model, ABMs are very attractive for modelers of various complex applications and disciplines (e.g., social sciences~\cite{bankes2002agent,macy2002factors}, ecology~\cite{grimm2005pattern,grimm2013individual}, archaeology~\cite{griffith2010hominids,graham2006networks}). However, 
the different existing formulations of ABMs appear to be quite inconsistent and lack rigorous theoretical foundation, 
%formulations of ABMs can be quite inconsistent and mostly lack a theoretical foundation,
which  motivates to develop standardized formats for these models~\cite{grimm2010odd,bankes2002agent,gilbert2002platforms}.   

In~\cite{ourepjpaper,conrad2018mathematical}, a mathematical formulation of an ABM is presented for systems of spatially distributed agents that move randomly in space and interact whenever they are close-by. The local interactions among agents trigger them to change their type, e.g., their opinion, or their infection status. The model  is formulated as a system of equations coupling the Markovian diffusion dynamics for the spatial agent movement  with Markov jump processes  for the type changes of each agent. Agent-based models  of this class can for instance be found as models for infection spreading~\cite{brockmann2013hidden,pastor2015epidemic},  innovation spreading~\cite{ourepjpaper,conrad2018mathematical} and chemical reactions~\cite{doi1976stochastic},  but also other classes of ABMs are commonly formulated as Markov processes \cite{izquierdo2009techniques}.
 
In general, these ABMs cannot be solved analytically, but require numerical methods for simulating realizations of the modeled dynamics. In~\cite{ourepjpaper}, 
a joint algorithm for the simulation of the spatial diffusion and the event-based interaction process is presented, building on the Euler-Maruyama scheme~\cite{higham2001algorithmic} and the  Temporal Gillespie algorithm~\cite{vestergaard2015temporal}. 
As the model formulation is stochastic,  many Monte-Carlo (MC) simulations are required in order to make adequate predictions for observables of the system.
However, for most real-world dynamics, such MC simulations are intractable due to an explosion of the costs for increasing agent numbers. 

\begin{figure}[htb!]
    \centering
    \includegraphics[width=\textwidth]{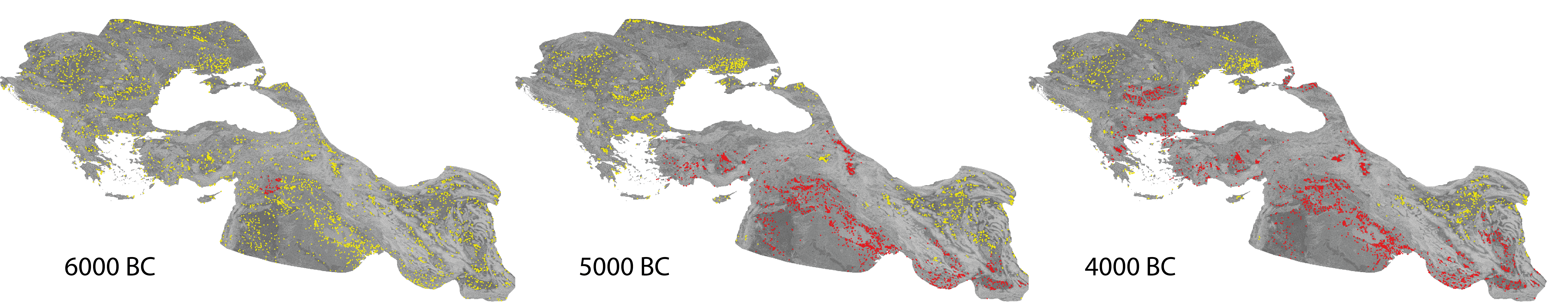}
    \caption{Snapshots of one realization of the wool-bearing sheep innovation spreading in ancient times~\cite{ourepjpaper,conrad2018mathematical}. Around $6000$ BC the innovation starts diffusing from the area of Tell Sabi Abyad in Near East towards Europe. Agents that do not know about the innovation are marked as yellow dots and red dots indicate agents that already have adopted the innovation.}
    \label{fig:sheep}
\end{figure}

In Figure \ref{fig:sheep}, we plot one realization of a MC simulation for innovation spreading in ancient times~\cite{ourepjpaper,conrad2018mathematical}. The underlying ABM model was built to study the spreading of the wool-bearing sheep in the Near East and Southeast Europe between $6000$ and $4000$ BC. This stochastic model considers $N=4000$ agents that diffuse in a complex, data-driven suitability landscape with non-trivial boundary conditions. Due to the high computational costs of simulating this model and the need for running many such MC simulations, developing new approaches for model reductions with a small approximation error is essential.

For spatial systems of many interacting agents, one can reduce the model complexity by replacing the micro-scale model of individual agents by a meso-scale model of stochastic agent densities, leading to a system of coupled stochastic partial differential equations (SPDEs)~\cite{dean1996langevin,kawasaki1973new,bhattacharjee2015fluctuating,kim2017stochastic}.
Classically, in the limit of infinitely many agents, stochastic effects become negligible and the dynamics can be described by a macro-scale model of deterministic reaction-diffusion PDEs~\cite{pearson1993complex,fisher1937wave}. In a real-world system, however, the number of agents is typically finite, such that the dynamics are still intrinsically random and stochastic modeling approaches are more appropriate in order to reflect the uncertainty and fully cover the emerging phenomena of the underlying ABM.  

In this paper, after reviewing and slightly generalizing  the stochastic ABM from~\cite{ourepjpaper,conrad2018mathematical}, we will formulate a system of SPDEs as a meso-scale approximation to the ABM, see~\cite{helfmann2018stochastic}.  
The model is an extension of the Dean-Kawasaki model~\cite{dean1996langevin,kawasaki1973new} that describes the transport in space of the stochastic agent densities and includes interactions between the agent densities for the different types~\cite{kim2017stochastic,bhattacharjee2015fluctuating}. The derived distribution-valued equations belong to the class of Dean-Kawasaki type problems, for which well-posedness is still unresolved due to several mathematical difficulties~\cite{cornalba2018regularised,lehmann2018dean,fehrmanwell}.

 Since in the end we are interested in the discretization and simulation of the model, we will, in order to deal with these issues, interpret the system of SPDEs in its weak form and consider ``a regularization by  discretization". 
 Henceforth, we propose an algorithm for the efficient sampling of trajectories of the coupled SPDE system by first discretizing in space using the Finite Element (FE) method, thereby approximating the unknown distributions by piece-wise polynomial functions. The resulting system of SDEs is well-posed and can be further discretized in time using the Euler-Maruyama method.
 
As  a meso-scale approach, simulation of the system of SPDEs is much faster than of the corresponding ABM, since the number of coupled equations is drastically reduced from scaling with the number of agents to scaling with the number of different agent types.
Further, the Finite Element method can easily deal with complicated domain  boundaries \cite{de2015finite}, which is crucial for many agent-based models and real-world situations (see Figure~\ref{fig:sheep}). Other approaches, e.g., using Finite Volume schemes have been intensively studied~\cite{bhattacharjee2015fluctuating,kim2017stochastic} and benefit from the property of agent number conservation.
 
 Last, we will study the ABM and the reduced SPDE model numerically on a toy example by comparing the computational effort and investigating the approximation quality of the reduced SPDE model with respect to the ABM.
The focus of the present  paper is on the mapping from agent-based models to SPDE models and their approximation quality and not on a numerical analysis of the approximation, still we will also shortly reflect on the  consistency of numerical solutions to the SPDE. 
 
\section{Agent-based Model}
Agent-based models  are micro-scale computational models describing the system of interest
in terms of discrete entities, called agents. An agent represents an autonomous entity such as a chemical particle, a person, a group of people or an organization. The individual behaviour of each agent is traced in space and time, including its  interactions with other agents and the environment.
The hope is that  from the interplay of the local behaviour of agents on the microscopic scale, global patterns emerge on a larger system scale ~\cite{helbing2012agent,macy2002factors,grimm2005pattern}. Emergence is likewise often stated as ``the whole is more than the sum of its parts" due to the interactions among the parts (i.e., the agents) resulting in global structures~\cite{bonabeau2002agent}.  

In this section, we will review and slightly generalize the agent-based model from ~\cite{ourepjpaper,conrad2018mathematical}, which is formulated as a system of  stochastic  processes describing  the position dynamics and the interaction rules for   agents. The position dynamics of agents are described by Brownian motion with a drift term  pushing them towards suitable parts of their local  surroundings. This is a useful description in the case that agents  move randomly and slightly jittery in space whilst following some (underlying) energy landscape, e.g., when describing the mobility of humans in the ancient world, the movement of animals,  or the diffusion of chemical particles~\cite{doi1976stochastic}. 

Further, agents can interact according to a set of predefined rules whenever they are close in space. These interactions cause agents to change their type influenced by the state of other near-by agents (e.g., their infection status~\cite{kermack1927contribution, pastor2015epidemic}, their current opinion~\cite{holme2006nonequilibrium}, whether they use an innovation~\cite{conrad2018mathematical,ourepjpaper}, whether they have a certain information). The type changes are commonly written as 
$T_1 + T_3 \rightarrow T_2+T_3$ 
(inspired by the notation for chemical reactions), where $T_1,T_2,T_3$ denote different agent types. 

Even though the model formulation is very simple,  it generates a vast range of different and complex outcomes on a larger system scale. 
The presented ABM is based on the Doi model ~\cite{doi1976stochastic} for reacting and locally diffusing chemical particles. 
It has similarities to   the SIR model for infection spreading~\cite{kermack1927contribution, pastor2015epidemic} and to Brownian agents~\cite{schweitzer2002modelling} and builds on the ABM for innovation spreading in the ancient world as introduced in~\cite{conrad2018mathematical,ourepjpaper} and Figure \ref{fig:sheep}. 

We will give a detailed description of  the agent-based model in Section \ref{sec:abm}, and explain a method for its simulation in Section \ref{sec:abm_sim}. 

\subsection{Model Formulation}\label{sec:abm}
We consider a system of $N$ agents, where each agent $i=1,\dots,  N$ is characterized by its type $Y_i(t)$ and its position $X_i(t)$ at time $t\in [0,T]$. 
The agents' positions can be continuous values in a given domain $D\subseteq \R^d$, whereas the type of an agent is a discrete feature denoted by values in $\{1,\dots, N_T\}$. Thus the state of the $i^\text{th}$ agent at time $t$ is given by
$
(X_i(t),\ Y_i(t)) \in D \times \{1,\dots,    N_T\}.
$
We write the state of the whole system of $N$ agents as
$
(X(t), Y(t)) = (X_i(t),\ Y_i(t))_{i=1}^N.$

In the model, we are following every agent $i=1,\dots,  N$ individually and track the evolution in time of its position state and type.  The position dynamics and the interactions between agents leading to type changes are described in the following sections.

\subsubsection{Position Dynamics of Agents}
 We assume that agents are able to move and change their position in the domain $D\subseteq \R^d$. Thereby agents are  taking  their local surroundings into account, in such a way, that they are attracted to near-by regions that are suitable for them and refrain from unsuitable parts of the domain.
 
We can straightforwardly model this by letting the agents follow the gradient of a potential landscape, the so-called suitability landscape. The suitability landscape $V$ indicates the attractivity of the environment and gives  an incentive for agents to prefer or avoid certain near-by parts of the domain.
 Valleys of the suitability landscape  correspond to attractive regions and peaks and divides correspond to unsuitable areas that are moreover difficult to overcome. The suitability landscape can be constructed on the basis of data and expert knowledge, see e.g.,~\cite{ourepjpaper,conrad2018mathematical}. 
Additionally, we include  randomness in the agents' motion to account for other unknown incentives for positional changes and to allow agents to be exploratory or make mistakes in their evaluation of the environment. 

Thus, the change of the position $X_i(t) \in D \subseteq \R^d$ of every agent $i=1 ,\dots,  N$ is governed by the Itô diffusion process 
\begin{equation}\label{eq:diffusion_process}
dX_i(t) = -\nabla V(X_i(t))dt + \sigma dB_i(t),
\end{equation}
with gradient operator $\nabla = \left(\frac{\partial}{\partial x_1},\dots, \frac{\partial}{\partial x_d}\right)^T,$ suitability landscape $V: D \subseteq \R^d \rightarrow \R$, diffusion constant\footnote{Note that within physical sciences the term \textit{diffusion constant} or \textit{diffusion coefficient} typically refers to  $\frac{\sigma^2}{2}$, however, in the theoretical context of stochastic differential equations these terms are also used for $\sigma$ itself~\cite{oksendal2013stochastic,evans2012introduction}.}  $\sigma\in\R$  and $
B_i(t)$ denoting independent standard Brownian motions in $\R^d$. We impose reflecting boundary conditions  in the case where $D$ is bounded and thereby ensure that the agents' positions are within $D$ for all times.  

\remark{The position dynamics of agents in many modeling scenarios are interdependent such that agents tend to group together in space and form clusters, while also keeping some distance from each other in order to avoid spatial overlap and crowding~\cite{conrad2018mathematical,ourepjpaper}. To account for this, another drift term has to be added to the position SDE (\ref{eq:diffusion_process}), s.t. the dynamics are governed by
\begin{equation}\label{eq:diffusion_process_attr_rep}
dX_i(t) = -\left(\nabla V(X_i(t))+ \nabla U_i(X(t))\right)dt + \sigma dB_i(t).
\end{equation}
Thereby each agent $i$ experiences an additional force derived from the attraction-repulsion potential~$U_i$
\begin{equation}\label{eq:att-rep-pot}
 U_i: D^N \subseteq \R^{d \times N}\rightarrow \R ,\ X(t) \mapsto \sum_{ j \neq i} u\left(\lVert X_i(t)-X_j(t) \rVert\right),
 \end{equation}
where we sum over all pair-wise attraction-repulsion potentials $u:\R_{\geq 0}\rightarrow \R$ 
between agent $i$ and $j \neq i$.
The pair-wise potentials $u$ are inspired by interatomic potentials from Physics (e.g., Lennard-Jones potential,  Buckingham potential). Attraction between pairs of agents occurs whenever agents at long distances are driven towards another, and repulsion appears when agents are forced apart at short distances. Agents are thus searching for an optimal balance between forming clusters of agents on the one hand and distributing in space on the other hand.}

\subsubsection{Interaction Rules}
Agents can change their type according to a set of $N_R$ interaction rules $\{R_r\}$, $r=1,\dots,  N_R$. The type changes happen at a certain rate and whenever an agent is in proximity of specific other agents that can influence  the agents' type. We consider rules $\{R_r\}$  that can be written as the type change\footnote{In principle, the interaction rules could be of a more complicated form, e.g., by including the death and birth of agents such as $T_s \rightarrow \emptyset$ and $ \emptyset \rightarrow T_s$, by allowing  both agents to change their type $T_s + T_{s''} \rightarrow T_{s'}+T_{s'''}$ or by letting a larger number of agent types take part in an interaction. A straightforward mapping of these more complicated interaction rules is also possible for the SPDE formulation that will be introduced in Section \ref{sec:model_SPDE},   we keep these simple interaction rules only for notational simplicity. } 
\begin{equation}\label{eq:rule}
R_r:\ T_s +  T_{s''} \rightarrow T_{s'}+T_{s''} \ \text{with} \ s,s',s'' \in \{1,\dots, N_T\}
\end{equation}
happening at the fixed influence rate $\gamma_r$  and triggered by a close-by agent of type $T_{s''}$.

For example when modeling the spreading of an innovation among humans ~\cite{ourepjpaper,conrad2018mathematical}, we can assume that every agent takes one of the two discrete innovation states: $T_1$ for a non-adopter or $T_2$ for an adopter of the innovation. By defining one simple interaction rule $R_1$: $T_1 + T_2 \rightarrow 2 \ T_{2}$, we can model that adopters  pass on the innovation to non-adopters at a fixed rate whenever they are in contact. 

Pairs of agents are in contact if they are within a distance $d_{\text{int}}$ of each other. Given the changing positions of agents (\ref{eq:diffusion_process}), we construct a time-evolving  network between agents (represented by the set of nodes) that are in contact (given by the edges) at time $t$. 
The network is fully determined by a time-dependent adjacency matrix $\mathcal{A}(t)=(\mathcal{A}_{ij}(t))_{i,j=1}^N$ with entries $\mathcal{A}_{ij}(t)=1$ if $i \neq j,\ \|X_i(t)-X_j(t)\| \leq d_{\text{int}}$, and 0 else.

As a next step, we are interested in describing the type change process $\{Y_i(t)\}_{t \in [0,T] }$ of each agent $i$. If agent $i$ at time $t$ is of type $T_s$, $s\in \{1,\dots,  N_T\}$, we denote this by $Y_i(t) = s$. Type changes for agent $i$ are modeled as Markov jump processes 
on $\{1,\dots,  N_T\}$ 
with time-dependent transition rates. The transition rates are changing in time since they depend both on the proximity of other agents and their types. The transition rate function $\lambda_i^r(t)$ gives the accumulated rate for agent $i$ of type $T_s$ to change its type according to interaction rule $R_r$ and is proportional to the constant influence rate $\gamma_r$ and to the number of neighbors of agent $i$ that trigger interaction $R_r$, i.e., the number of agents of type $T_{s''}$ in rule (\ref{eq:rule}).

Consequently, we write the transition rate function for  agent $i$ at time $t$ following interaction rule $R_r:\ T_s +  T_{s''} \rightarrow T_{s'}+T_{s''}$ as 
\begin{equation}\label{eq:inf_rate}
  \lambda_i^r(t)  = \lambda_i^r(\mathcal{A}(t),Y(t)) = \gamma_r  \ \sum_{j=1}^N \mathcal{A}_{ij}(t) \  \1_{\{s'' \}} (Y_j(t)) \ \1_{\{s \}} (Y_i(t)), 
\end{equation}
where $\1_{B}$ is the indicator function defined as
$\1_{B}(x)= 1$ if $x \in B$ and 0 else.
Finally, the type change process $\{Y_i(t)\}_{t \in [0,T]}$ of agent $i$   can  be expressed as
\begin{equation}\label{eq:typechanges}
Y_i(t)=Y_i(0) + \sum_{r=1}^{N_R} \mathcal{P}_i^r \left( \int_0^t \lambda_i^r(t') dt' \right) v_r,
\end{equation}
where $Y_i(0)$ denotes the initial type of agent $i$, $\mathcal{P}_i^r$ denote i.i.d. unit-rate Poisson processes and the type change vector is denoted by $v=(v_r)_{r=1}^{N_R}$ ($v_r=s'-s$ for $R_r$ as given in \eqref{eq:rule}  above).

\subsubsection{Formulation of the Agent System Dynamics}
Putting together Equations (\ref{eq:diffusion_process}) and (\ref{eq:typechanges}), we can describe the evolution in time of the agent states by the following coupled equations
\begin{align}\label{eq:coupled_process}
X_i(t) &=X_i(0) - \int_0^t \nabla V(X_i(t'))dt' + \sigma \int_0^t dB_i(t') \nonumber \\
Y_i(t)&=Y_i(0) + \sum_{r=1}^{N_R} \mathcal{P}_i^r \left( \int_0^t \lambda_i^r(t') dt' \right) v_r
\end{align}
 for agents $i=1,\dots,  N$.

Since we cannot solve the coupled equations (\ref{eq:coupled_process}) analytically, we will in the following explain how to accurately discretize and efficiently simulate trajectories of the dynamics.
\subsection{Simulation}\label{sec:abm_sim}
For each agent $i$, the type change process depends on the positions and on the types of all other agents via the time-dependent transition rate functions $\lambda_i^r(t)$.
If we additionally include attraction-repulsion forces between agents, then further the motion of all agents is intrinsically intertwined.
For the design of the simulation algorithm we have to take this into account and therefore propose to simulate both the position  and type changes \eqref{eq:coupled_process} in parallel. 

The most straightforward approach is to discretize time $[0,T] $ with sufficiently small $\Delta t$ time steps. The Euler-Maruyama method can be employed to accurately and efficiently discretize the position  SDE~\cite{kloeden,ourepjpaper}, whereas when discretizing the Markov jump processes we replace rates by probabilities in each time step~\cite{conrad2018mathematical,ourepjpaper}.  
In contrast to time-continuous simulation schemes for the interactions, like the standard Gillespie algorithm, such a time discretization can evoke interaction conflicts in each integration step, namely when two type change events are chosen to occur for the same agent. In order to circumvent such conflicts, we define for each type $T_s$ the index set
\begin{equation}\label{Rs}
    \mathcal{R}_s := \{r \in \{1,...,N_R\}: R_r:\ T_s +  T_{s''} \rightarrow T_{s'}+T_{s''} \} 
\end{equation}
of rules that change this type, and stipulate that each agent of type $T_s$ can undergo maximally one event $R_r$, $r\in \mathcal{R}_s$, per time step.
The resulting simulation approach can be summarized as follows.

For each time step $t_k=k \Delta t$, $k=0,\dots,   K-1$:
\begin{itemize}
\item[1.] For each agent $i=1,\dots, N$: the positions are advanced according to $$
X_i(t_{k+1}) =X_i(t_k) -   \nabla V(X_i(t_k)) \Delta t + \sigma\sqrt{\Delta t} \ \zeta_{i,k} 
$$
with i.i.d. $\zeta_{i,k} \sim \mathcal{N}(0,1)$.
\item[2.] For each $i=1,\dots,N$, set $s:=Y_i(t_k)$ (the type of agent $i$ at time $t_k$) and calculate $$\Lambda_i(t_k) := \sum_{r \in \mathcal{R}_{s}}  \lambda_i^r(t_k)$$ with $\mathcal{R}_{s}$ defined in \eqref{Rs}.  Then, a type change of agent $i$ occurs with probability $$1-\exp(-\Lambda_i(t_k)\Delta t).$$ In case the agent changes its type, decide subsequently by which of the rules, choosing $R_r$ with probability $$\lambda_i^r(t_k)/\Lambda_i(t_k)$$ for each $r\in \mathcal{R}_{s}$.
\end{itemize}
A few remarks are in order: 
Note that this procedure does not restrict the total number of interaction events per agent and time step to one. An agent  may still undergo several interaction events in parallel as long as maximally one of these events changes its type (i.e., for an agent of type  $T_{s''}$ several events of the form \eqref{eq:rule} are possible). Further, since the underlying interaction dynamics are time-continuous, a more accurate simulation~\cite{fennell2016limitations} of the type changes would simulate them  continuously in time using a Gillespie algorithm~\cite{vestergaard2015temporal}, for the coupled simulation approach in this case we refer the reader to~\cite{ourepjpaper}.
 
Regarding the computational cost, we note that  the algorithm scales more than linearly with the number of agents $N$,  since in each time step, we need to compute the contact network $\mathcal{A}(t_k)$, i.e., the pair-wise distances between agents. The brute-force approach of computing all pair-wise distances scales like $\mathcal{O}(N^2)$ for $N$ agents, whereas clever algorithms can lower the cost  to $\mathcal{O}(N^{5/3})$ ~\cite{verlet1967computer}.
These distances are also required when including attraction-repulsion forces, since the pair-wise attraction-repulsion potential $u\left(\lVert X_i(t)-X_j(t) \rVert\right)$
depends on the pair-wise distances between agents. 

The proposed simulation algorithm produces a single   trajectory of the joint stochastic process $\{X(t),Y(t)\}_{t\in [0,T]}$. In order to obtain reliable statistics of the dynamics and to estimate the process' distribution and its moments, several thousand repeated simulations are required. 
In case of a system with large numbers of agents this 
causes problems since the simulation cost  scales badly with increasing agent numbers $N$. Parameter studies are numerically not tractable in this case. Many real-world systems actually contain very large numbers of agents, which motivates to consider a less complex modeling approach  in terms of stochastic PDEs, and thereby offers a  cheap method for simulating trajectories.

\section{Model Reduction to a System of Stochastic PDEs}
\begin{figure}[htb!]
    \centering
    \includegraphics[width=0.6\textwidth]{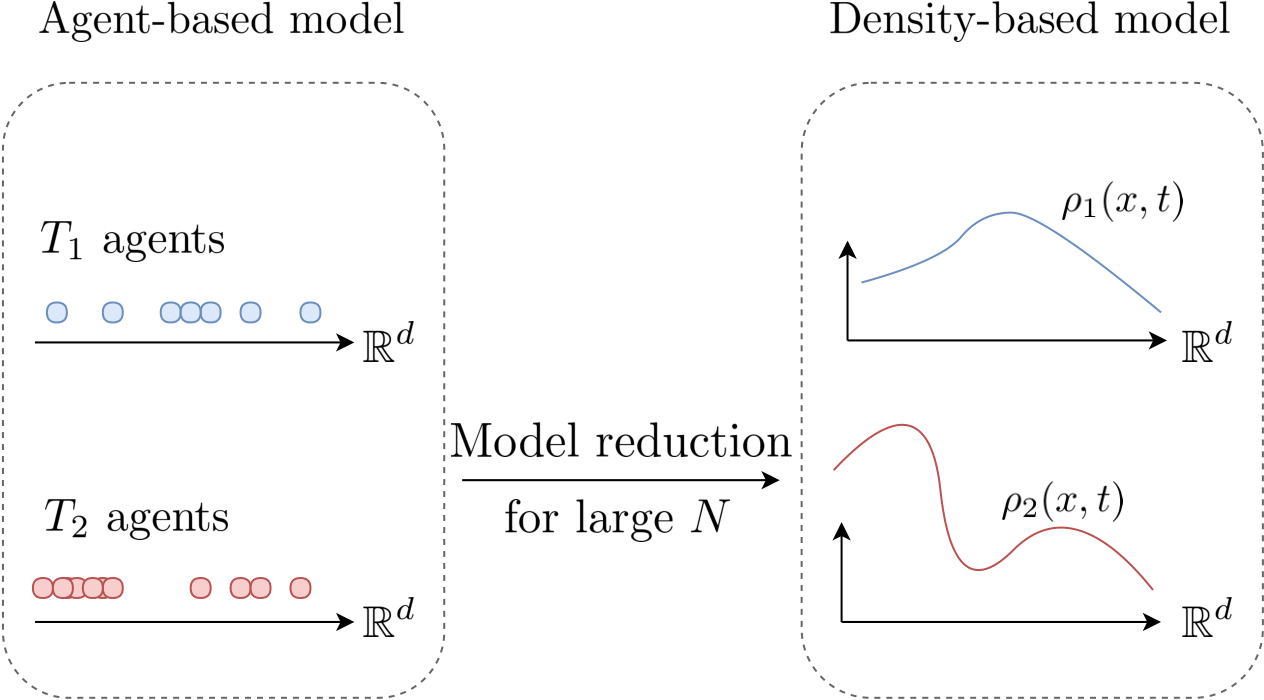}
    \caption{Model reduction from the agent-based formulation in terms of $N$ discrete agents to the SPDE description in terms of agent number densities $\rho_s(x,t)$ for types $s=1,\dots,N_T$.}
    \label{fig:approxmodel}
\end{figure}

Instead of studying the evolution of every individual agent, in the following we study the transport in space and the interactions between agents in terms of number densities  for each type, such that the model complexity and simulation costs are drastically reduced.
This reduced model, given in the form of a system of stochastic PDEs, is a combination of the model for  position changes proposed by Dean and Kawasaki~\cite{dean1996langevin,kawasaki1973new} (also termed Dean-Kawasaki model), and  added terms for the interactions  between agent densities~\cite{bhattacharjee2015fluctuating,kim2017stochastic}. It approximates the ABM for systems of many, but finitely many agents.
Stochasticity still emerges from the systems' inherent randomness due to the finite number of agents.  But agents become indistinguishable among their type and we lose the individual agent labels, see Figure \ref{fig:approxmodel}.

Applications of similar SPDE models  include chemical pattern formation~\cite{bhattacharjee2015fluctuating,kim2017stochastic} and models for bacterial populations~\cite{grafke2017spatiotemporal}, but to our knowledge they have not yet been applied for modeling systems of humans. 

We will start by proposing and explaining the features of the reduced SPDE model  in Section \ref{sec:model_SPDE}. In Section \ref{sec:discr}, we will tackle its efficient numerical discretization by means of the Finite Element method. 

\subsection{Model Formulation} \label{sec:model_SPDE}
We are considering a model describing the stochastic evolution of agent number densities (or number concentrations) for each agent type $s=1,\dots,  N_T$. The (stochastic) agent number density 
\[
\rho_s:  D \times [0,T] \mapsto \R_{\geq 0}
\]
 is defined on the domain $D\subseteq \R^d$, the time interval $[0,T]$ and with probability space $(\Omega, \mathcal{F}, \mathbb{P})$\footnote{Actually the number density should be $\rho_s(\omega, x,t)$ but we do not write the explicit dependence on $\omega \in \Omega$. The agent number density is random since it solves a stochastic PDE.}. Integrating the number density $\rho_s$ over the domain yields the number of agents of type $T_s$, which we denote by $N_s$.
 
The densities $(\rho_1(x,t),\dots,  \rho_{N_T}( x,t))$ evolve due to diffusion and drift in the suitability landscape  and because of the set of interaction rules. The  temporal changes of  
$
\rho_s( x,t)$ for $s=1,\dots,   N_T$ 
are given by the following stochastic partial differential equation (SPDE)  interpreted in the  It\^{o} sense
\begin{equation}\label{eq:spde_op}
\frac{\partial \rho_s( x,t)}{\partial t} = \mathcal{D} \rho_s( x,t) + \mathcal{I} \rho_s( x,t) 
\end{equation}
with stochastic diffusion operator $\mathcal{D}$ (see Section \ref{sec:diff_op}) and stochastic interaction operator $\mathcal{I}$ (Section \ref{sec:int_op}). 
For a fixed sample $\omega \in \Omega$, the agent density is a realization of a stochastic process solving the SPDE.

\subsubsection{Diffusion of the Agent Densities}\label{sec:diff_op}
The diffusion operator in Equation (\ref{eq:spde_op}) is given by~\cite{dean1996langevin,kawasaki1973new, kim2017stochastic,bhattacharjee2015fluctuating}
\begin{equation}
\mathcal{D} \rho_s( x,t) :=  \frac{\sigma^2}{2} \Delta \rho_s(x,t) +  \nabla \cdot (\nabla V(x) \rho_s(x,t)) +\sigma \nabla \cdot \left(\sqrt{\rho_s(x,t)} Z_s^{\mathcal{D}}(x,t) \right)   \label{eq:diffusion_spde}
\end{equation}
with diffusion constant $\sigma \in \R$, Laplace operator ${\Delta = \sum_{l=1}^d \frac{\partial^2}{\partial x_l^2}}$ and suitability landscape ${V: D \subseteq \R^d \mapsto \R}$. $Z_s^{\mathcal{D}}(x,t)=\left(Z_{s,1}^{\mathcal{D}}(x,t),\dots,  Z_{s,d}^{\mathcal{D}}(x,t)\right)$ 
denotes a $d$-dimensional vector of space-time white noise (STWN) for the diffusion, i.e., a mean zero  process that is uncorrelated in space and time $$\E \left( Z_{s,j}^{\mathcal{D}}(x,t) Z_{s',j'}^{\mathcal{D}}(x',t') \right)= \delta_{jj'} \delta_{ss'} \delta(x-x') \delta(t-t'),$$
where $\delta_{ij}$ denotes the Kronecker Delta and $\delta(x-y)$ denotes the Dirac Delta distribution.

The diffusive part of the SPDE evolves a number density of many agents and is responsible for  the diffusive transport in space with drift in the suitability landscape $V(x)$. 
Its deterministic part $ \partial_t \rho_s  =\frac{\sigma^2}{2} \Delta \rho_s +  \nabla \cdot (\nabla V \rho_s)$ reminds of the Fokker-Planck equation that describes the evolution of a density of infinitely many non-interacting particles that are  diffusing with drift in a potential.
Here though, we consider a large but finite number of diffusing particles, as such the system is still intrinsically random and the derivation leads to an additional noise term  fluctuating around zero. The derivation~\cite{dean1996langevin,kawasaki1973new} of $\frac{\partial \rho_s}{\partial t} = \mathcal{D} \rho_s$ from the position SDE \eqref{eq:diffusion_process} is based on the following idea: We define the empirical agent density for a single agent $i$ as  
\begin{equation}\label{eq:dirac}
\rho^i(x,t) =  \delta(x-X_i(t)),
\end{equation}
where $\delta(x-X_i(t))$ denotes a Dirac Delta distribution placed at position $X_i(t)$. Then, using It\^{o}'s Formula, the SDE  (\ref{eq:diffusion_process}) for agent $i$ is transformed into an SPDE describing the temporal evolution of $\rho^i(x,t)$. 
We sum the SPDEs over $i=1,\dots,N_s$ to get an equation for the density of all agents  of type $T_s$, i.e., for
$\rho_s(x,t)=\sum_{i=1}^{N_s} \rho^i(x,t).$ 
Aiming at a closed-form equation for $\rho_s(x,t)$ that does not depend on $\rho^i(x,t)$ explicitly, Dean~\cite{dean1996langevin}  in his formal computations approximates the obtained conditional Gaussian noise term by a Gaussian noise of the same mean and covariance that just depends on $\rho_s(x,t)$. 
 
This multiplicative noise term $\sigma \nabla \cdot (\sqrt{\rho_s} Z_s^{\mathcal{D}}) $  is non-linear  and comes with several mathematical problems, making the question of existence and uniqueness of the solution still an open problem~\cite{cornalba2018regularised,lehmann2018dean}. First, the noise term is given in a divergence form. In order to interpret this term, we will consider the weak formulation and use partial integration to transfer the divergence operator to the test functions. Second,  the meaning of the term $\sqrt{\rho_s}$ remains undefined. Namely, the unknown $\rho_s$ is defined as a sum of Dirac Delta distributions, hence the meaning of its square root is unclear. We will deal with this problem by approximating the distribution $\rho_s$ by a function from a Finite Element space, for which the square root is well defined. 

Different notions of solutions to (\ref{eq:diffusion_spde}) that would be natural to consider exist, such as a martingale solution~\cite{lehmann2018dean} or a path-wise kinetic solution~\cite{fehrmanwell}. In order  to apply the Finite Element method and to reformulate the divergence of the noise term, in our setting it is more appropriate to consider the weak formulation (in the PDE sense) that will be derived in \ref{sec:weak}.

\remark{
We can extend the diffusion operator to include attraction forces between pairs of agents at long ranges and repulsion forces at short ranges, similar as in Equation (\ref{eq:diffusion_process_attr_rep}). Denoting the pair-wise attraction-repulsion potential between two agents at positions $x$ and $y$ by $u(\lVert x-y \rVert)$, the SPDE is extended by one term (the derivation can be found in~\cite{dean1996langevin})
\begin{align} \label{eq:attraction_repulsion}
\mathcal{D} \rho_s( x,t) =& \frac{\sigma^2}{2} \Delta \rho_s(x,t) +  \nabla \cdot (\nabla V(x) \rho_s(x,t))  \nonumber \\ &+ \nabla \cdot \left(\rho_s(x,t) \int_D \left(\sum_{s'=1}^{N_T}\rho_{s'}(y,t)\right) \nabla u(\lVert x-y \rVert) dy\right) +\sigma \nabla \cdot \left(\sqrt{\rho_s(x,t)} Z_s^{\mathcal{D}}(x,t)\right).  
\end{align}
The additional term models the diffusion of the agent density $\rho_s(x,t)$ in the aggregated attraction-repulsion potential of the density of all agents $\sum_{s'}\rho_{s'}(x,t)$. 
Including this term, the diffusing densities are all coupled to each other.

For analytical simplicity, we will not further consider attraction and repulsion forces in the remaining paper. But it should be possible to extend the discretization of the SPDE straightforwardly.}

\subsubsection{Interactions between Agent Densities}\label{sec:int_op}
The interaction operator $\mathcal{I}$ of the system of SPDEs (\ref{eq:spde_op}) accounts for the local transport  of some density  between the different agent densities $\rho_s(x,t)$, $s=1,\dots,N_T$ in accordance with the set of   interaction rules $R_r,\ r=1,\dots,N_R$~\cite{kim2017stochastic,bhattacharjee2015fluctuating}
\begin{equation}
\mathcal{I} \rho_s(x,t) := \sum_{r=1}^{N_R} \nu_s^r \left( a^r(\rho(x,t)) + \sqrt{a^r(\rho(x,t))} Z_r^{\mathcal{I}}(x,t)\right).  \label{eq:interaction_spde}
\end{equation} 
Thereby the coefficient $\nu_s^r$  describes the discrete number change of the type $T_s$ agents involved in the $r^{\text{th}}$ interaction rule\footnote{In the chemistry literature, $\nu_s^r$ is called the stoichiometric coefficient of  type $T_s$ due to the $r^{\text{th}}$ reaction. Further, for rules $R_r$ of the form \eqref{eq:rule}, $\nu_s^r$ and $v_r$ (the type change coefficient from the ABM) are related via $v_r=\sum_s \nu_s^r s$.}, $a^r(\rho(x,t))$ is the transition rate function for densities according to rule $R_r$, and 
space-time white noise for the $r^\text{th}$ interaction is denoted by $Z_r^{\mathcal{I}}(x,t)$ with covariance 
$$\E \left(Z_r^{\mathcal{I}}(x,t) Z_{r'}^{\mathcal{I}}(x',t') \right) = \delta_{rr'} \delta(x-x') \delta(t-t').
$$ 
Similar as in the ABM, the rate function takes into account the local amount of the two types of agents taking part in the interaction.  For interaction rule $R_r:\ T_s + T_{s''} \rightarrow T_{s'}+T_{s''}$ with $s,s',s'' \in \{1,\dots, N_T\}$, the transition rate function at time $t$ for $T_s$ agents reads
\begin{equation}\label{eq:tr_SPDE}
a^r(\rho(x,t)) 
= \rho_{s}(x,t) \gamma_r  \int_{B_\text{int}(x)}   \rho_{s''}(x',t) dx'
\end{equation}
since all the mass of type $T_s$ agents  that are placed at $x\in D$ can interact at rate $\gamma_r $ with all agents of  type $T_{s''}$ that are   within the closed ball $B_{\text{int}}(x)$ of radius $d_\text{int}$~\cite{donev2018efficient,erban2009stochastic}. The number of  $T_{s''}$ agents within that radius is simply given by the integral of the number density over the interaction neighbourhood.

To get a better understanding of these coefficients and functions, let us return to our innovation spreading example: agents of type $T_1$ and $T_2$ are interacting according to the rule $R_1:\ T_1 + T_2 \rightarrow 2 \  T_2$. Since for each interaction the number of type $T_1$ agents decreases by one agent and the number of type $T_2$ agents increases by one agent, we have $\nu_1^1 =-1,\ \nu_2^1=1$. The transition rate function~\eqref{eq:tr_SPDE} for the interaction between two agent densities is proportional to the density $\rho_1(x,t)$ of $T_1$ agents at $x$, the number of $T_2$ agents within the interaction radius around $x$, and the rate $\gamma_1 $  such that 
$a^1(\rho(x,t))=\gamma_1  \ \rho_1(x,t) \int_{B_\text{int}(x)}   \rho_{2}(x',t) dx'$
in this example. 

In order to derive the interaction operator of the SPDE (\ref{eq:interaction_spde}), the Poisson random variable for interactions in the agent-based description (\ref{eq:typechanges}) has been replaced by a Gaussian random variable with the same mean and variance~\cite{kim2017stochastic}. This approximation is only valid for large 
$a^r(\rho(x,t))$, i.e., especially for large populations, and is closely related to the
approximation that has been done for well-mixed systems of interacting species leading to the Chemical Langevin equation (CLE)~\cite{gillespie2000chemical} in the context of the so-called large population limit. Equation (\ref{eq:interaction_spde}) includes spatial information and can be viewed as a spatial extension of the CLE.

\remark{In~\cite{kim2017stochastic}, in the context of numerical simulation of the model, the Gaussian noise is again replaced by Poisson noise in each grid cell. The reason for this is that the Gaussian approximation is only valid in the large population limit, i.e., for large numbers of agents in each grid cell. Switching back to Poisson noise  helps decreasing the resulting lack of accuracy and preventing negative values for densities in numerical simulations; however, the procedure is rather ad-hoc and no precise understanding of its effect is available.  }

\subsubsection{Complete Dynamics}\label{sec:complete}
For the system of SPDEs to be fully determined, boundary  and initial conditions have to be specified. The domain boundary should be Lipschitz continuous. Since in our  model agents cannot leave the domain, we require no-flux boundary conditions on $\delta D$. The initial data $\rho_{s,0}(x)$ has to be non-negative and such that $\int_D \rho_{s,0}(x) dx = N_s$.
Then, for agent types $s=1,\dots,  N_T$ the system of  SPDEs 
reads
\begin{align} \label{eq:spde}
\frac{\partial \rho_s( x,t)}{\partial t} &= \mathcal{D} \rho_s( x,t) + \mathcal{I} \rho_s( x,t) & \text{on}\ D \times [0,T] \nonumber\\
\nabla \rho_s( x,t) \cdot \hat{u}(x) &= 0 & \text{on}\ \delta D \times [0,T] \nonumber\\
\rho_s(x,0) &= \rho_{s,0}(x) & \text{on}\ D \times \{0\},
\end{align} 
where $\hat{u}$  denotes the normal to the boundary $\delta D$.

Note that the choice of the boundary conditions can be done differently. In particular, besides asking that the deterministic flux is zero across the boundary, one can also impose a zero stochastic flux at the boundary, as it is done in \cite{de2015finite}. In the weak formulation of the equation (see Section \ref{sec:weak}) this choice of boundary conditions would result in all boundary integral terms to vanish. However, this would not significantly change  our computational results, since we choose 
the potential $V$ such that it keeps the agent number density away from the boundary. For this reason and the problem of mathematically strictly defining deterministic and stochastic flux boundary conditions, we impose only a  deterministic zero-flux boundary condition.

As we have now   described the system of SPDEs, there are some further points worth noting. We can observe that the diffusion operator $\mathcal{D}$ from Equation (\ref{eq:diffusion_spde}) conserves the number of agents of each type because of the divergence operator form and assuming no-flux boundary conditions.  It is solely responsible for the transport of the density in space. The interaction part (\ref{eq:interaction_spde}) of the SPDE shifts the agent number density locally between the different agent types, but conserves the overall number of agents of all types, i.e., $\sum_{s=1}^{N_T} \int_D \rho_s(x,t) dx=N$ is conserved since  $\sum_s \nu_s^r = 0$ for each interaction rule of the form  \eqref{eq:rule}. In the case of more complicated interaction rules (e.g., the birth and death of agents), this is usually not the case.

We can also have a closer look at how the random forcings scale for increasing agent numbers, i.e., larger values of the density $\rho_s(x,t)$. The noise terms are scaled by a square-root factor of the agent density, whereas all the other (deterministic) terms are scaled by the density. Thus for the number of agents approaching zero, the noise dominates the SPDE.
For the number of agents going to infinity, the noise terms become unimportant and could be neglected. Then, the SPDE could be replaced by a PDE model.  Our numerical examples in Section \ref{sec:num} will display this behaviour.

The analysis of the well-posedness and existence of solutions to this SPDE system is not investigated enough~\cite{fehrmanwell,lehmann2018dean,cornalba2018regularised}. In this paper though, we are concerned with a discretization of the SPDE system. We will in the following show how to formally derive the weak formulation of the system of SPDEs (Section \ref{sec:weak}) forming the basis for discretizing the SPDEs with the Finite Element method (Section \ref{sec:discr}). 

\subsubsection{Weak Formulation}\label{sec:weak}
The SPDE  (\ref{eq:spde}) is properly interpreted as an  integral equation in time. 
For this, we will introduce the Cylindrical Wiener process $\{W(\cdot,t)\}_{t\in[0,T]}$ ~\cite{lord2014introduction,da2014stochastic} as
\begin{equation} \label{eq:wiener}
W(x,t)=\sum_{m=1}^\infty \chi_m (x) B_m(t),\end{equation}
where $\{\chi_m(x)\}_{m\in\mathbb{N}}$ is any orthonormal basis of $L^2(D)$ and $B_m(t)$ are i.i.d. Brownian motions  in $\mathbb{R}$. The Cylindrical Wiener process is a stochastic process that is  Brownian in time and white (i.e., uncorrelated) in space, such that its time derivative turns out to be space-time white noise $Z(x,t)=\partial_t W(x,t)$. By introducing the Cylindrical Wiener process expansion \eqref{eq:wiener} into the SPDE (\ref{eq:spde}), we  can  use the integral theory for stochastic processes in space and time~\cite{lord2014introduction,da2014stochastic}. Additionally, this offers a straightforward possibility for simulating realizations of $Z(x,t)$ by  truncating the expansion \eqref{eq:wiener} to a finite number of terms, numerically differentiating in time, and sampling i.i.d. Brownian motions $B_m(t)$.
We denote by $W_s^{\mathcal{D}}$ and $W_r^{\mathcal{I}}$ the noise terms that correspond to diffusive and interactive part respectively and we assume that they are independent.

In order to find solutions of the system of SPDEs, we will consider the weak solution framework (in the PDE sense). For the derivation of the weak form, we interpret \eqref{eq:spde} with the introduction of \eqref{eq:wiener} in the time integral sense,  multiply it by test functions $w(x)$, integrate over the domain $D$, and use partial integration. 
As usual, by making use of partial integration in deriving the weak form, we reduce the regularity requirements of the solution $\rho_s(x,t)$. Moreover, utilizing the no-flux boundary conditions $\nabla \rho_s(x,t) \cdot \hat{u}(x) =0$ on $\delta D$, we get
\begin{align}%\label{eq:green}
&\langle \Delta \rho_s(\cdot,t), w \rangle =   -\int_D \nabla \rho_s(x,t) \cdot \nabla w(x) \ dx + \int_{\delta D} (\nabla \rho_s(x,t) \cdot \hat{u}(x)) w(x) \ dx  = - \langle \nabla \rho_s(\cdot,t) , \nabla w \rangle. \nonumber
\end{align}
Here and in the following, we use the inner product notation $\langle u,v \rangle = \int_D u(x) v(x) \ dx$.

Further, we   shift the divergence operator from the space-time white noise onto the test functions. By using partial integration, we find
\begin{equation}\label{eq:stoch_weak}
\left\langle \nabla \cdot \left(\sqrt{\rho_s(\cdot,t)} dW_s^{\mathcal{D}}(\cdot,t)\right)  , w\right\rangle = - \left\langle \sqrt{\rho_s(\cdot,t)} dW_s^{\mathcal{D}}(\cdot,t)  ,\nabla  w  \right\rangle + \int_{\delta D} \left(\sqrt{\rho_s(x,t)} dW_s^{\mathcal{D}}(x,t) \cdot \hat{u}(x) \right) \ w(x)\ dx.
\end{equation}
With that, the weak formulation of (\ref{eq:spde}) consists of finding $\rho_s(x,t)$  
for all agent types $s=1,\dots,N_T$, such that
\begin{align}\label{eq:weak_3}
&\langle \rho_s(\cdot,t),w  \rangle =  \langle \rho_{s,0} , w  \rangle + \int_0^t \left( - \frac{\sigma^2}{2} \left\langle \nabla \rho_s(\cdot,t'), \nabla w  \right\rangle + \left\langle \nabla \cdot \left( \nabla V \rho_s(\cdot,t') \right), w  \right\rangle \right) dt' \nonumber\\ 
&+ \sum_{r=1}^{N_R} \nu_s^r \int_0^t \left\langle a^r(\rho(\cdot,t')) ,w  \right\rangle  dt'
-  \sigma\int_0^t \left\langle   \sqrt{\rho_s(\cdot,t')} dW_s^{\mathcal{D}}(\cdot,t')  ,\nabla  w  \right\rangle \nonumber\\  
&+\sigma \int_0^t \left( \int_{\delta D} \left(\sqrt{\rho_s(x,t')} dW_s^{\mathcal{D}}(x,t')  \cdot \hat{u}(x)\right) w(x)\ dx \right)
+   \sum_{r=1}^{N_R} \nu_s^r \int_0^t  \left\langle \sqrt{a^r(\rho(\cdot,t'))} dW_r^{\mathcal{I}}(\cdot,t'),w  \right\rangle
\end{align}
holds for all $ w(x)$ and for all $t\in [0,T]$.

The integrals with respect to the Cylindrical Wiener process have to be understood as follows
$$\int_0^t \left\langle   \sqrt{\rho_s(\cdot,t')} dW_s^{\mathcal{D}}(\cdot,t')  ,\nabla  w  \right\rangle
=\sum_{l=1}^d \sum_{m=1}^\infty \int_0^t \left\langle \sqrt{\rho_s(\cdot,t')} \chi_m  , \frac{\partial w  }{\partial x_l} \right\rangle dB_{s,m,l}^\mathcal{D}(t')
$$
$$\int_0^t \int_{\delta D} \left(\sqrt{\rho_s(x,t')} dW_s^{\mathcal{D}}(x,t')  \cdot \hat{u}(x)\right) w(x)\ dx 
=\sum_{l=1}^d \sum_{m=1}^\infty \int_0^t \left(  \int_{\delta D} \sqrt{\rho_s(x,t')} \chi_m(x) \hat{u}_l(x) w(x) dx \right) dB_{s,m,l}^\mathcal{D}(t')
$$
$$\int_0^t  \left\langle \sqrt{a^r(\rho(\cdot,t'))} dW_r^{\mathcal{I}}(\cdot,t'),w  \right\rangle
= \sum_{m=1}^\infty \int_0^t \left\langle a^r(\rho(\cdot,t')) \chi_m , w   \right\rangle dB_{r,m}^\mathcal{I}(t').
$$
 
This formal derivation of the  weak form \eqref{eq:weak_3} is the foundation for discretizing the system of SPDEs in the following section. 

\subsection{Discretization of the System of SPDEs}\label{sec:discr}
Previously, we introduced a model (\ref{eq:spde}) for the evolution of stochastic agent densities $\rho_s: D \times [0,T]  \rightarrow \mathbb{R}_{\geq 0} $ that are being transported in space and interacting with each other. 
The existence and uniqueness of solutions  to the SPDE is still an open question, especially in the last few years much research has been focused on this~\cite{fehrmanwell,lehmann2018dean,cornalba2018regularised}.
 
Here we instead propose a Finite Element discretization of the weak formulation of the system of SPDEs  (\ref{eq:weak_3}) via a finite spatial ansatz space and truncation of the noise expansion. The result of this Galerkin discretization will be the (finite) system of SDEs~\eqref{eq:galerkin5}. For such SDEs rigorous statements concerning the existence and uniqueness of solutions exist: if the drift and noise intensity terms satisfy appropriate conditions regarding Lipschitz continuity and growth at infinity, then the solution exists and is unique~\cite[Chapter 2]{Mao2008} and moreover stable with respect to perturbations of drift and noise intensity~\cite[Chapter 4]{Mao2008}.  

The resulting system of SDEs can be further discretized in time, s.t. we arrive at (a sequence of) matrix equations that have to be solved at each time step.

The general steps of the spatial discretization are as follows~\cite{lord2014introduction}:
\begin{enumerate}
    \item We project the solutions onto the finite-dimensional space $\tilde{V}$ (Section \ref{sec:space}).
    \item We truncate the  expansion of the Cylindrical Wiener process $W^M(x,t)=\sum_{m=1}^M \chi_m (x) B_m(t)$ to $M$ terms   (Section \ref{sec:truncation}).
    \end{enumerate}
    
In order to arrive at a fully discrete system we will perform a third step at last:
    \begin{enumerate}\setcounter{enumi}{2}
    \item We discretize in time using the Euler-Maruyama scheme (Section \ref{sec:time}).
\end{enumerate}
The discretization of the system of SPDEs will allow us to very efficiently simulate trajectories that approximate the dynamics of the agent-based model.  

Instead of discretizing the SPDE model with a Finite Element approach, one can also employ the Finite Volume method~\cite{kim2017stochastic,donev2010accuracy,delong2013temporal}.
Here, we use the Finite Element method since it   has the advantage that one can in principle treat very complicated  domains and boundaries \cite{de2015finite} which is needed for many real-world models.  Further, due to the weak form interpretation, we can make sense
of the divergence operator acting on the space-time white noise.

\subsubsection{Space Discretization}\label{sec:space}
We consider a simplicial discretization (e.g., by triangles, rectangles) $\tilde{D}=\tilde{D}(h)  $  of the (possibly curvilinear) computational domain $D \subseteq \mathbb{R}^d$.  Here $h$ is defined as $h:=\max_{E \in \mathcal{T}} h_E$, where $h_E$ is the diameter of the simplex $E$ and $\mathcal{T}$ is the union of all simplices. We let $\tilde{V}=\tilde{V}(h)$ be a finite-dimensional element space consisting of continuous piece-wise linear functions defined on $\tilde{D}$ and spanned by its nodal basis (e.g., hat functions) $\{\phi_i:\tilde{D}\rightarrow \R\}_{i=0}^n$. In the following, the inner product on the discretized domain $\tilde{D}$ is given by $\langle u,v \rangle_h= \int_{\tilde{D}} u(x) v(x) \ dx$ and the normal to the discretized boundary is denoted by $\tilde{u}$.

The Finite Element  method  then reduces the problem (\ref{eq:weak_3}) to finding $ \tilde{\rho}_s(\cdot,t)\in \tilde{V}$ for each agent type $s=1,\dots,N_T$ such that 
\begin{align}\label{eq:spde44}
&\left\langle  \tilde{\rho}_s(\cdot,t),\phi_i  \right\rangle_h = \langle \tilde{\rho}_{s,0} , \phi_i  \rangle_h + \int_0^t \left( - \frac{\sigma^2}{2} \left\langle \nabla  \tilde{\rho}_s(\cdot,t'), \nabla \phi_i  \right\rangle_h + \left\langle \nabla \cdot \left( \nabla V  \tilde{\rho}_s(\cdot,t') \right), \phi_i  \right\rangle_h \right) dt' \nonumber\\ 
&+ \sum_{r=1}^{N_R} \nu_s^r \int_0^t \left\langle a^r(\tilde{\rho}(\cdot,t')) ,\phi_i  \right\rangle_h  dt'
-  \sigma\int_0^t \left\langle   \sqrt{ \tilde{\rho}_s(\cdot,t')} dW_s^{\mathcal{D}}(\cdot,t')  ,\nabla  \phi_i(x) \right\rangle_h \nonumber\\ 
&+\sigma \int_0^t \left( \int_{\delta \tilde{D}} \left(\sqrt{ \tilde{\rho}_s(x,t')} dW_s^{\mathcal{D}}(x,t')  \cdot \tilde{u}(x) \right) \phi_i(x)\ dx \right)  
+   \sum_{r=1}^{N_R} \nu_s^r \int_0^t \left\langle \sqrt{a^r(\tilde{\rho}(\cdot,t'))} dW_r^{\mathcal{I}}(\cdot,t'),\phi_i  \right\rangle_h \nonumber \\
&\forall t\in[0,T] \ \text{and for all test functions } \phi_i\in \tilde{V}, i=0,\dots,n.
\end{align}
To get the initial data of the discretization, we have to project the given   $\rho_{s,0}(x)$ onto $\tilde{V}$,  i.e.,
$$\tilde{\rho}_{s,0}(x) = \sum_{j=0}^n \langle \rho_{s,0} , \phi_j  \rangle \phi_j(x).$$
With  $\{\tilde{\rho}_s(\cdot,t)\}_{t\in[0,T]}$  being a $\tilde{V}$-valued stochastic process, we can next expand a realization of $\{\tilde{\rho}_s(\cdot,t)\}_{t\in[0,T]}$  as a linear combination of the basis functions $\{\phi_j\}_{j=0}^n$  with time-dependent (and random) coefficients $\beta_{s,j}(t)$,  i.e.,
\begin{equation}\label{eq:real_expan}
 \tilde{\rho}_s(x,t) = \sum_{j=0}^n \beta_{s,j}(t) \phi_j(x).\end{equation}
We define matrices 
\begin{equation}\label{eq:matrix_C}
C := \left(C_{ij}\right)_{i,j=0}^n = \left( \langle \phi_j ,\phi_i  \rangle_h \right)_{i,j=0}^n ,
\end{equation}
\begin{equation}\label{eq:matrix_A}
A := \left( A_{ij} \right)_{i,j=0}^n = \left( \frac{\sigma^2}{2} \left\langle \nabla \phi_j , \nabla \phi_i  \right\rangle_h - \left\langle \nabla \cdot \left( \nabla V \phi_j  \right), \phi_i  \right\rangle_h \right)_{i,j=0}^n,\end{equation}
the coefficient vector 
\begin{equation}\label{eq:beta}
\beta_s(t):=(\beta_{s,j}(t))_{j=0}^n,
\end{equation}
and the vector of the deterministic (non-linear) interaction term coupling the densities via $a^r(\tilde{\rho}(x,t))$ 
\begin{equation}\label{eq:F}
F_s(t)=(F_{s,i}(t))_{i=0}^n:=\left(  \sum _{r=1}^{N_r} \nu^r_{s} \left\langle a^r(\tilde{\rho}(\cdot,t)) ,\phi_i  \right\rangle_h \right)_{i=0}^n.
\end{equation}
Inserting the expansion (\ref{eq:real_expan}) into (\ref{eq:spde44}) and using the defined quantities (\ref{eq:matrix_C}) - (\ref{eq:F}),
we finally arrive at
\begin{align}\label{eq:space_disc}
&\sum_{j=0}^n C_{ij} d\beta_{s,j}(t) = - \sum_{j=0}^n A_{ij} \beta_{s,j}(t) \ dt
+ F_{s,i}(t) \ dt 
-  \sigma   \left\langle \sqrt{ \tilde{\rho}_s(\cdot,t)} dW_s^{\mathcal{D}}(\cdot,t)  ,\nabla  \phi_i  \right\rangle_h \nonumber\\
&+\sigma \int_{\delta \tilde{D}} \left(\sqrt{ \tilde{\rho}_s(x,t)} dW_s^{\mathcal{D}}(x,t)  \cdot \tilde{u} (x)\right) \phi_i(x)\ dx 
+  \sum_{r=1}^{N_R} \nu_s^r  \left\langle \sqrt{a^r(\tilde{\rho}(\cdot,t))} dW_r^{\mathcal{I}}(\cdot,t),\phi_i  \right\rangle_h, \ \nonumber \\
&\forall t\in[0,T], \ \forall  i=0,\dots,   n,
\end{align}
which is still understood as a time integral equation in the sense of (\ref{eq:spde44}).

\subsubsection{Truncation of the Noise Expansion} \label{sec:truncation}
By truncating the noise expansion \eqref{eq:wiener} to $M$ dimensions, we project the Cylindrical Wiener processes onto the finite-dimensional space spanned by  $\{\chi_m\}_{m=1}^M$. The choice of the truncation threshold $M$ remains an open question (e.g.,~\cite[Theorem 10.41]{lord2014introduction}), in our numerical experiments we will follow \cite{lord2014introduction} and  choose $M=n$.

The truncated expansion for the interaction noise reads 
\begin{equation}\label{eq:noise1}
W^{\mathcal{I},M}_r(x,t):=\sum_{m=1}^M \chi_m (x) B^{\mathcal{I}}_{r,m}(t),
\end{equation}
where $B_{r,m}^{\mathcal{I}}(t)$  denote i.i.d. Brownian motions in $\mathbb{R}$. Analogously for the diffusion noise, we replace $W_s^\mathcal{D}(x,t)$ by the truncated noise expansion\footnote{The choice of the threshold $M$ could in general be different for \eqref{eq:noise1} and \eqref{eq:noise2}.}
\begin{equation}\label{eq:noise2}
W^{\mathcal{D},M}_s(x,t):=\sum_{m=1}^M \chi_m (x) B^{\mathcal{D}}_{s,m}(t),
\end{equation}
where $B_{s,m}^{\mathcal{D}}(t)=(B_{s,m,l}^{\mathcal{D}}(t))_{l=1,...,d}$  denote i.i.d. Brownian motions in $\mathbb{R}^d$.
 
Then by defining the vectors
$$B^{\mathcal{I},M}_r(t):=(B^{\mathcal{I}}_{r,m}(t))_{m=1}^{M}$$
$$B_{s,l}^{\mathcal{D},M}(t):= \left( B_{s,m,l}^{\mathcal{D}}(t) \right)_{m=1}^M$$
and matrices
\begin{equation}\label{eq:GD}
    G_{s,l}^{\mathcal{D}}(t)_{im} := - \sigma   \left\langle \sqrt{\tilde{\rho}_s(\cdot,t)}  \chi_m ,\frac{\partial\phi_i}{\partial x_l}  \right\rangle_h +  \sigma \int_{\delta \tilde{D}} \left(\sqrt{ \tilde{\rho}_s(x,t)} \chi_m(x) \ \tilde{u}_{l} (x)\right) \phi_i(x)\ dx \ \text{for } l=1,\dots,d
\end{equation}
and 
\begin{equation}\label{eq:GI}
G_{s,r}^{\mathcal{I}}(t)_{im} := \nu^r_{s}  \left\langle \sqrt{a^r(\tilde{\rho}(\cdot,t))}  \chi_m ,\phi_i  \right\rangle_h, \ \text{where} \ i=0,\dots,n,\ m=1,\dots,M,
\end{equation}
our Galerkin approximation finally reads
\begin{align}\label{eq:galerkin5}
C d \beta_s(t) &=  (- A \beta_s(t)
+ F_s(t)) \ dt  + \sum_{l=1}^d G_{s,l}^{\mathcal{D}}(t) \ dB_{s,l}^{\mathcal{D},M} (t) + \sum_{r=1}^{N_R} G_{s,r}^{\mathcal{I}}(t) \ dB_r^{\mathcal{I},M} (t). 
\end{align}
Next, we will discretize this system of SDEs \eqref{eq:galerkin5} in time using the (backward or semi-implicit) Euler-Maruyama method.

\subsubsection{Time Discretization}\label{sec:time}
The last step is to discretize in time, we divide the time interval $[0,T]$ into $K$ intervals of sufficiently small size $\Delta t$. Denoting functions at time $t_k= k \Delta t$ by a subscript $k$, e.g., $\beta_s(t_k)=\beta_{s,k}$, the semi-implicit Euler-Maruyama time-discretization (implicit in the linear terms, explicit in the non-linear coupling term) of (\ref{eq:galerkin5}) is the recursion for $k=0,\dots,   K-1$
\begin{align} \label{eq:sim_spde}
 \beta_{s,k+1} &= (C + A \Delta t)^{-1} \left(C \beta_{s,k} + F_{s,k}  \Delta t + \sum_{l=1}^d G_{s,l,k}^{\mathcal{D}}\Delta B_{s,l,k}^{\mathcal{D}} + \sum_{r=1}^{N_R} G_{s,r,k}^{\mathcal{I}}\Delta B_{r,k}^{\mathcal{I}}  \right).
\end{align}
The Brownian increments $$
\Delta B_k  =  \left(\int_{t_k}^{t_{k+1}} dB_m(t) \right)_{m=1}^M =  \left( B_m(t_{k+1})-B_m(t_k)\right)_{m=1}^M =  \left( \sqrt{\Delta t} \zeta_{m,k}\right)_{m=1}^M
$$ 
have to be sampled for each time step $t_k$ and for each $s=1,\dots,N_T,$ $l=1,\dots,d$ and  $r=1,\dots,N_R$   by drawing i.i.d. $\zeta_{m,k} \sim \mathcal{N}(0,1)$.

It is known that this time discretization converges globally~\cite{Mao2012} with convergence order $1/2$ for $\Delta t\to 0$. 

\remark{
One numerical problem is that the agent density can become negative in the simulations~\cite{kim2017stochastic}. Since we want to be able to take square roots of the density, we can tackle this problem in the implementation by working with $\text{max}\{0, \beta_{s,j,k}\}$ instead of $\beta_{s,j,k}$  in the expressions $F_{s,k}$, $G_{s,l,k}^{\mathcal{D}}$ and $G_{s,r,k}^{\mathcal{I}}$ of the discretization \eqref{eq:sim_spde}.
This and other reformulations of the SDE in the vicinity of $\beta_{s,j,k}=0$ are discussed broadly in the literature, e.g., in the context of chemical Langevin equations \cite{gillespie2000chemical}, see \cite{cCLE} for a rigorous approach to negative densities. 

Another undesired property of the Finite Element discretization is that the number of agents is not necessarily conserved, instead the best solution within the finite-dimensional solution space is found. Since we want the number of agents to be conserved, we can circumvent that problem by re-normalizing the density after each time step thus enforcing the total number of agents to remain constant. For conserving the number of agents, the Finite Volume method is a natural candidate method \cite{kim2017stochastic}. }
\section{Numerical Studies}\label{sec:num}
Many models of real-world dynamics pose challenges regarding their simulation due to a complex model formulation. The agent-based model as introduced in Section \ref{sec:abm} can be considered as the ground-truth model for a system of diffusing and interacting agents. 
However due to its high model complexity for large numbers of agents, simulations of the discretized ABM are often too expensive. The ABM is valid on all population scales but we expect its computational feasibility and necessity only on the smallest population scale.  
Reduced model descriptions are therefore needed, these reduced models should have a small approximation error not only when the system is in stationarity but also when it is out-of-equilibrium, as well as  being computationally much more efficient.  In Section \ref{sec:model_SPDE}, we presented a reduced model in terms of stochastic PDEs approximating the original ground-truth system with full spatial resolution and for large populations, whilst still including stochasticity.

In this section, we will illustrate both modeling approaches on a toy example of innovation spreading among human agents that has real-world applications  ~\cite{conrad2018mathematical,ourepjpaper}. After introducing the dynamics and studying each of the models in some detail for  a fixed number of  agents in Section \ref{sec:ex_abm} and \ref{sec:ex_spde},  and analysing the numerical discretization of the SPDE in \ref{sec:num_con},  we will in Section \ref{sec:ex_comp} compare the two modeling approaches for different population sizes regarding their computational cost and the approximation quality of the system of SPDEs to the ABM.

\subsection{Illustrative Example: Modeling with Agents}\label{sec:ex_abm}
\begin{figure}[htb!]
    \centering
    \includegraphics[width=0.5\linewidth]{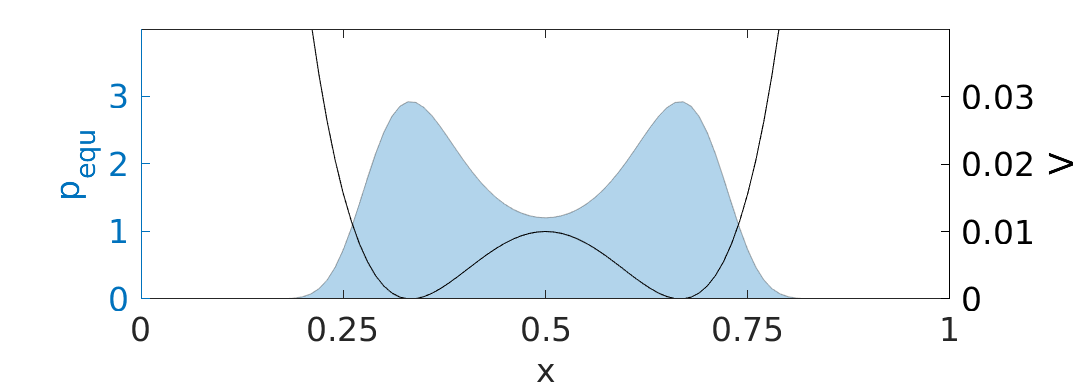}
    \caption{Suitability landscape $V$ and the corresponding unique equilibrium distribution  for the position of one agent $p_{\text{equ}}(x) =Z^{-1}\ \text{exp}\left(- 2 \sigma^{-2}  V(x) \right)$, the  distribution is normalized by $Z$.}
    \label{fig:eq_V}
\end{figure}

We study $N=1000$ agents diffusing  with diffusion constant $\sigma = 0.15$ in the double well landscape $V(x)=0.01\,(3.6(x-0.5)^2 - 0.1)^2$ on $D=[0,1]$, see Figure \ref{fig:eq_V}.
The suitability landscape $V$ is characterized by two minima centered at $x=\frac{1}{3}$ and $x=\frac{2}{3}$, corresponding to the most preferred spots for agents and a barrier between the two wells. Thus the agents will spend most of their time near the centers of the two wells while rarely transitioning between them.

On top of the position dynamics, agents are interacting whenever they are close-by. We consider the spreading of an innovation given by the rule $R_1: T_1 + T_2 \rightarrow 2\ T_2$ at fixed influence rate $\gamma_1=0.1$ \  whenever two agents of type $T_1$ (non-adopter of the innovation) and $T_2$ (adopter of the innovation) are within radius $d_\text{int}=0.002$ of each other.
Due to the small interaction radius, agents from one well can only interact with agents of the same potential well.

\begin{figure}[htb!]
\begin{subfigure}[c]{0.33\textwidth}
\includegraphics[ width=\linewidth]{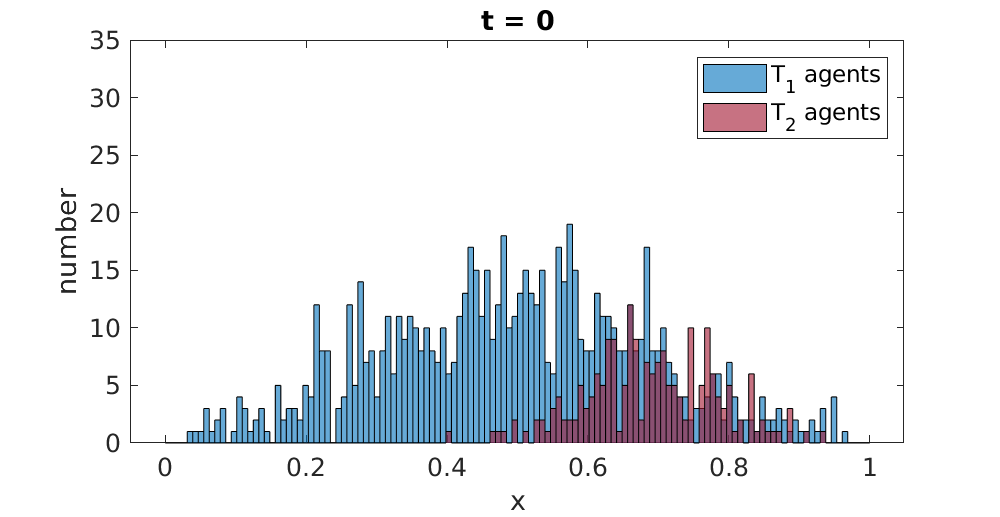}
\end{subfigure}
\begin{subfigure}[c]{0.33\textwidth}
\includegraphics[width=\linewidth]{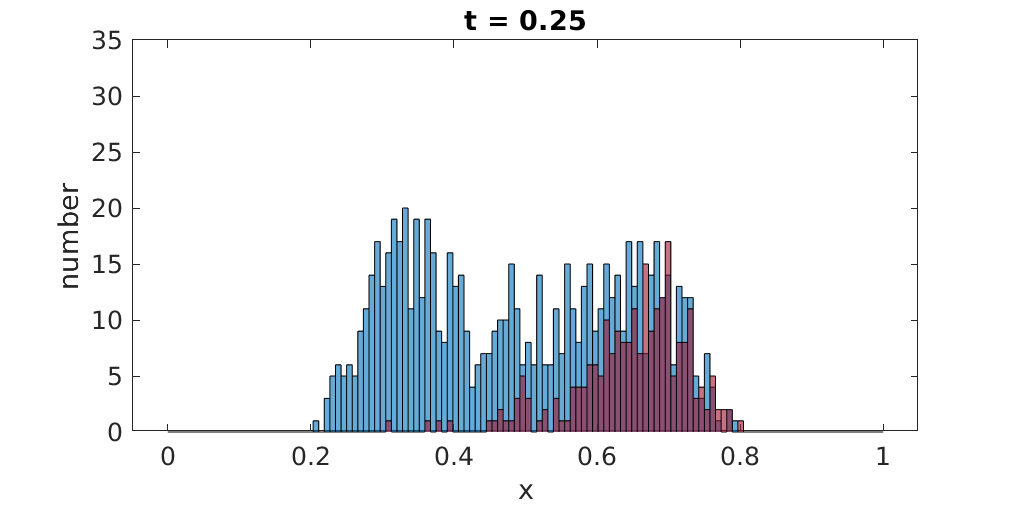}
\end{subfigure}
\begin{subfigure}[c]{0.33\textwidth}
\includegraphics[width=\linewidth]{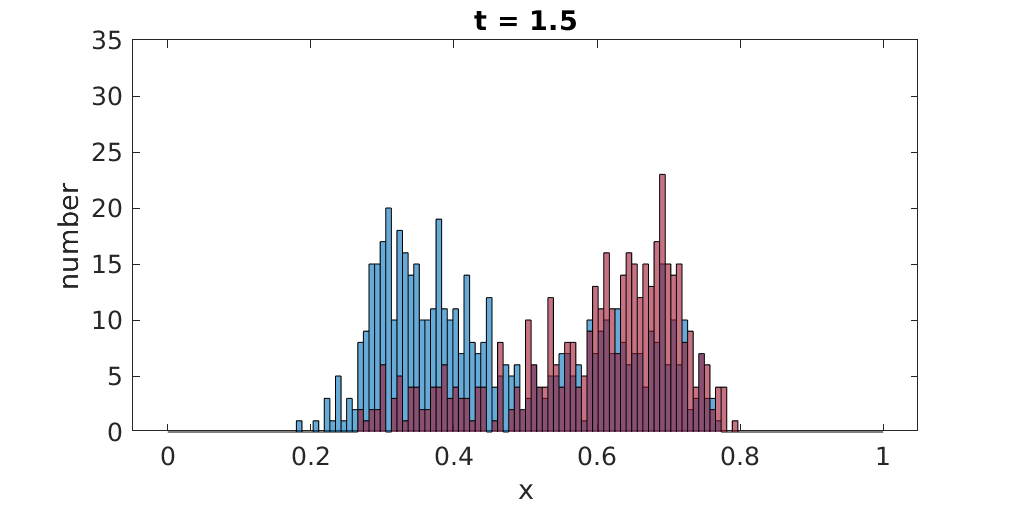}
\end{subfigure}

\begin{subfigure}[c]{0.33\textwidth}
\includegraphics[width=\linewidth]{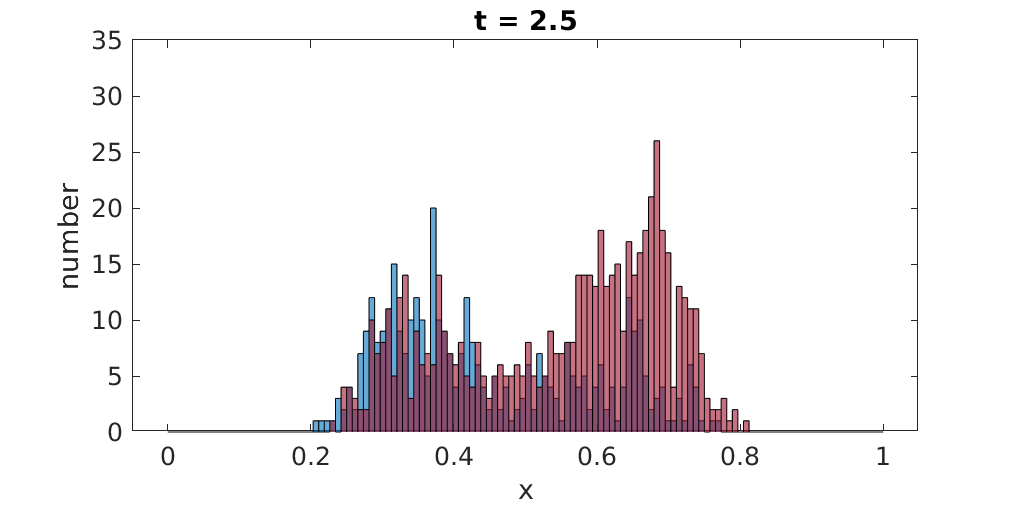}
\end{subfigure}
\begin{subfigure}[c]{0.33\textwidth}
\includegraphics[width=\linewidth]{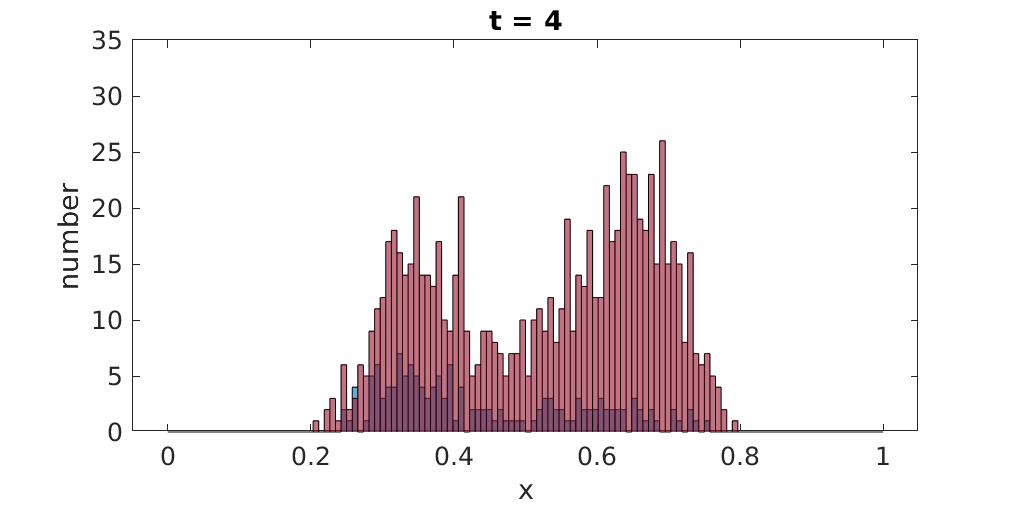}
\end{subfigure}
\begin{subfigure}[c]{0.33\textwidth}
\includegraphics[width=\linewidth]{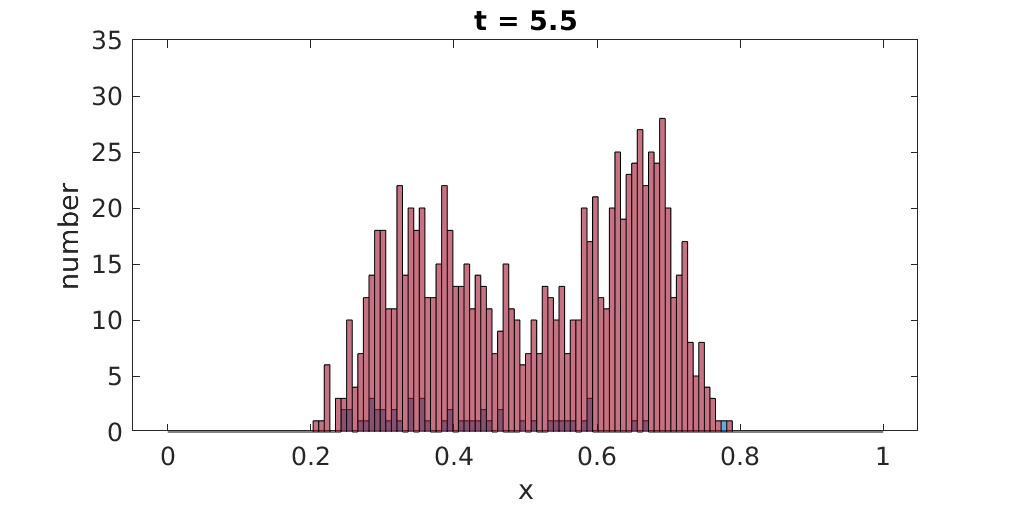}
\end{subfigure}
\caption{ A single realization of the spreading process in a double well landscape on $D=[0,1]$ with $1000$ agents and modeled by the ABM. The empirical density  of $T_1$ agents at time $t$ is given by the blue histogram while the density of  type $T_2$ agents is shown using the red histogram.}\label{fig:ex_abm}
\end{figure}

We will describe the dynamics by the coupled ABM equations \eqref{eq:coupled_process} and simulate discretized realizations of the process  with step size $\Delta t = 0.01$ as explained in Section \ref{sec:abm_sim}. 
We assume that at $t=0$, the positions of the initially $800$ agents of type $T_1$ are   normally  distributed  with mean $0.5$ and standard deviation $0.2$. Similarly, for the $200$ agents of type $T_2$, we draw the initial positions from the normal distribution $\mathcal{N}(0.7,0.1)$. 
Snapshots of one simulation of the spreading process are given in Figure \ref{fig:ex_abm}. 

We observe the following dynamics: At time $t=0$, the innovation starts spreading in the well centered at $x=\frac{2}{3}$. The agents quickly distribute near the attractive centers of the two wells. It takes some time until the innovation reaches the other well centered at $x=\frac{1}{3}$. But as soon as an adopter agent crosses the barrier for the first time, the innovation quickly gets adopted by all agents in the other well. At the final time  $t=5.5$, nearly all agents are of type $T_2$ and distributed according to the equilibrium distribution of the landscape $V$ (see Figure \ref{fig:eq_V}).

\subsection{Illustrative Example: Modeling with Agent Densities}\label{sec:ex_spde}
As in the previous section, we consider agents in a double well landscape that are passing on an innovation. But this time,
we will model the dynamics using the system of SPDEs that evolves the two agent densities $\rho_1(x,t)$ (density of non-adopters) and $\rho_2(x,t)$ (density of adopters).

To get the initial conditions, two approaches  are possible. Either the initial densities are given already in the form of a density, or in the case that the positions and types of the $N$ agents at time $t=0$ are given, we can construct the initial densities $\rho_{s,0}(x)$ by summing unit masses (e.g., normalized  narrow Gaussian functions, hat functions) placed at the position of each agent of type $T_s$. 
In this example, we will take the first approach by using (similar as in the previous model setting) that the agent positions are normally distributed, i.e.,  $\rho_{1,0}(x)$ and $\rho_{2,0}(x)$ are Gaussian functions integrating to the number of agents of each type ($N_1=800$ and $N_2=200$ respectively) at time $t=0$. The remaining parameters are chosen similar to the agent-based model.

Given the model setting, we can turn to the specification of the discretization. Numerically sampling trajectories of the SPDE system can be done by iteratively solving \eqref{eq:sim_spde} in parallel for both agent types $s=1,2$.
For the space discretization, we split the domain $D=[0,1]$ into equidistant grid cells of size $h=\frac{1}{n}$ (s.t. $\tilde{D}(h)=D)$,  supporting the hat functions $\phi_i(x)$, $i=0,\dots,n$, defined as
\[ \phi_i(x) := \begin{cases} 
      0 & x < x_{i-1} \text{ or }  x \geq x_{i+1} \\
      \frac{x-x_{i-1}}{h} & x_{i-1}\leq x < x_i \\
      1- \frac{x-x_i}{h} & x_i \leq x < x_{i+1}, 
   \end{cases}
\]
where $x_i=i h$. The hat functions are spanning our Finite Element space $\tilde{V}$. 
For the expansion of the cylindrical Wiener process, we choose the trigonometric system   {  $\{\sqrt{2} \sin(\pi m x)\}_{m=1}^{M}$ as an orthonormal basis for $L^2[0,1]$ which is zero at the domain boundary. By this choice, the stochastic boundary terms in the weak formulation vanish (see Equation \eqref{eq:stoch_weak}).} 
The matrices \eqref{eq:matrix_C}, \eqref{eq:matrix_A}, \eqref{eq:GD}, \eqref{eq:GI} can be assembled by analytical or numerical integration, or by defining a reference element and transforming from any grid cell to the reference grid cell in order to do the computations on the reference cell before transforming back~\cite{larsson2008partial}. 
Since each hat function is mostly zero throughout $\tilde{D}$, the resulting matrices will be sparse.
In the numerical time stepping, grid densities can become negative, therefore we work with $\text{max}\{0, \beta_{s,j,k}\}$ instead of $\beta_{s,j,k}$ whenever a square-root is taken and in the interaction rate $a^r(\tilde{\rho}(x,t))$. Further, the number of agents is not exactly conserved, which we circumvent  by re-normalizing the total density  to $N$ after each time step 
(see also the remark in Section \ref{sec:time}).

\begin{figure}[htb!]
\begin{subfigure}[c]{0.33\textwidth}
\includegraphics[width=\linewidth]{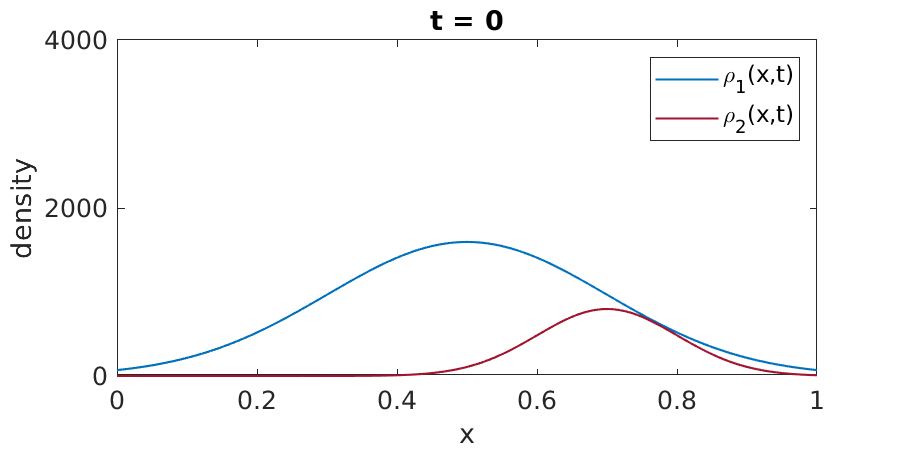}
\end{subfigure}
\begin{subfigure}[c]{0.33\textwidth}
\includegraphics[width=\linewidth]{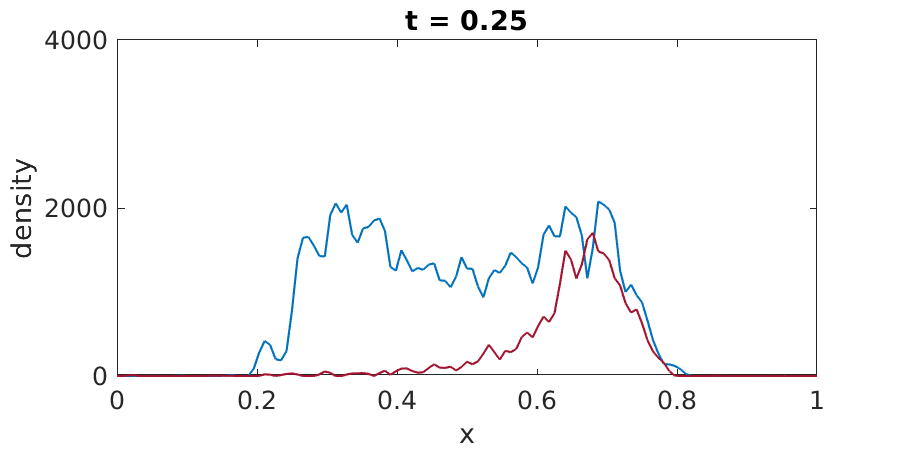}
\end{subfigure}
\begin{subfigure}[c]{0.33\textwidth}
\includegraphics[width=\linewidth]{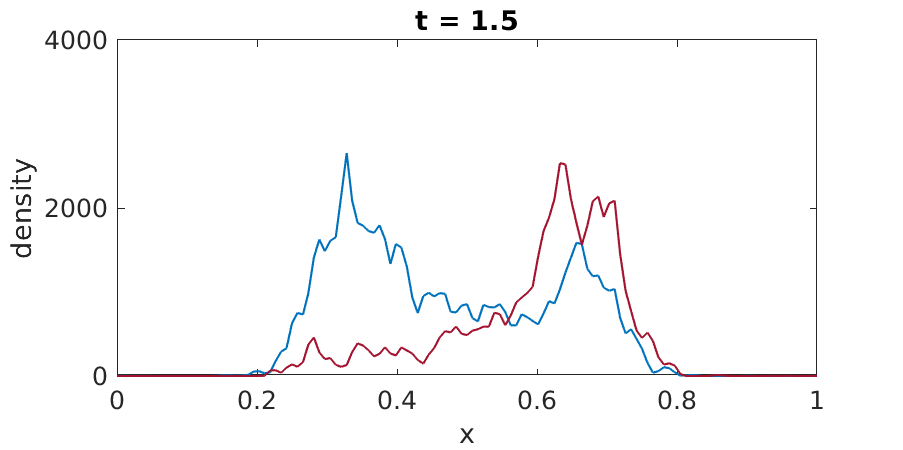}
\end{subfigure}

\begin{subfigure}[c]{0.33\textwidth}
\includegraphics[width=\linewidth]{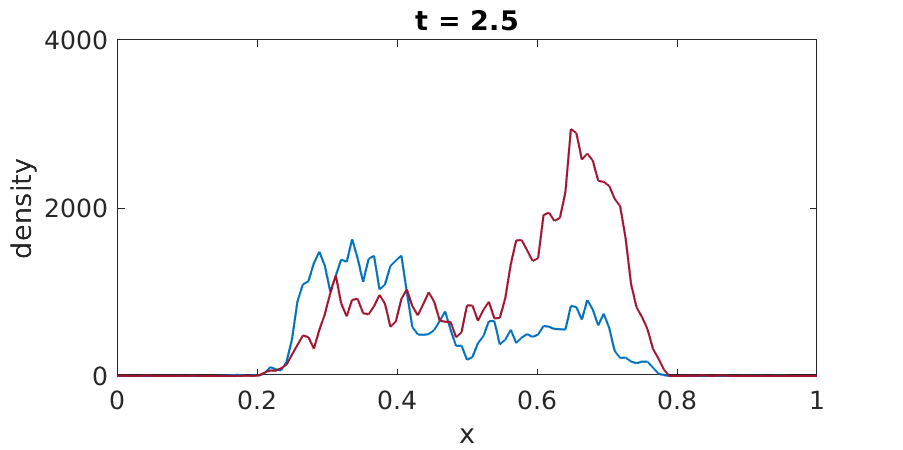}
\end{subfigure}
\begin{subfigure}[c]{0.33\textwidth}
\includegraphics[width=\linewidth]{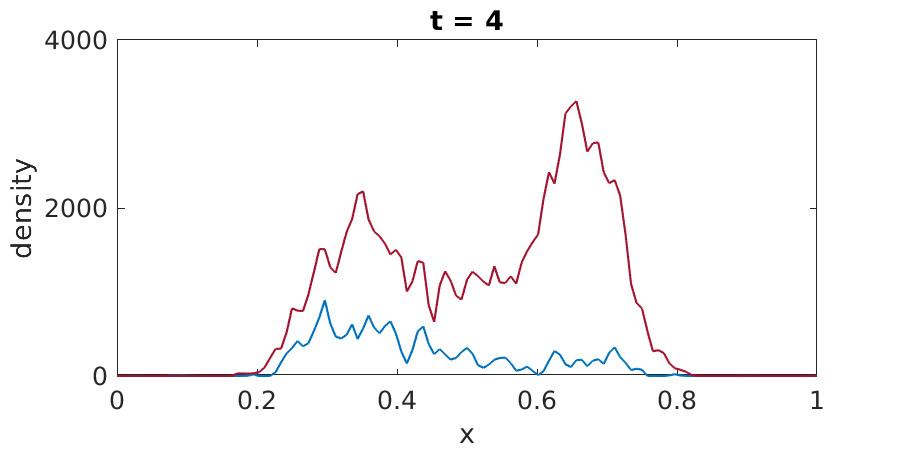}
\end{subfigure}
\begin{subfigure}[c]{0.33\textwidth}
\includegraphics[width=\linewidth]{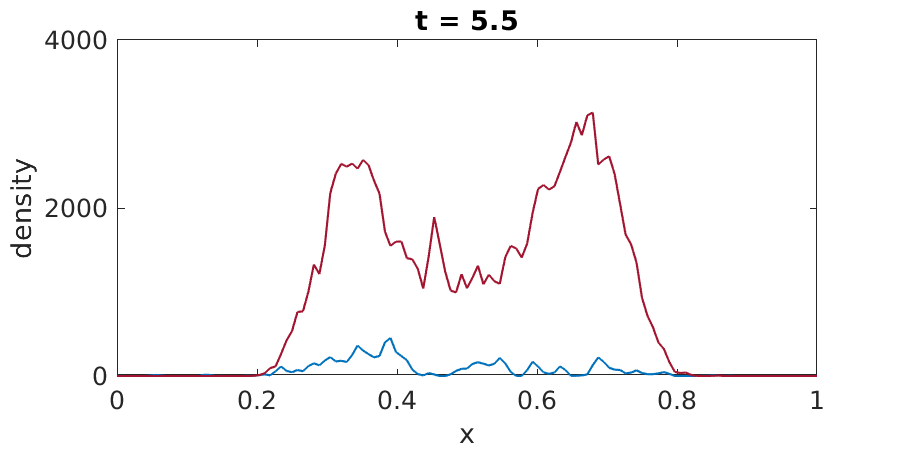}
\end{subfigure}
\caption{Realization of the innovation spreading process in a double well landscape  with $1000$ agents (SPDE approach). We plot the number density of non-adopters $\rho_1(x,t)$ and the density of adopters $\rho_2(x,t)$  at a few time instances. }\label{fig:ex_spde}
\end{figure}

Snapshots of one realization of the discretized SPDE solution are shown in Figure \ref{fig:ex_spde},  where as discretization parameters we choose  $n=M=128$ and time steps of size  $\Delta t=0.01$.

The emerging dynamics agree qualitatively with the ABM dynamics in Figure \ref{fig:ex_abm}. Again the agent densities are clustered around the two wells. The spreading of the innovation among   agents inside the same well is fast, but the spreading of the innovation from one well to another takes much longer.
We observe that the stochasticity inherent in the model is still visible on the global scale, the agent densities are seemingly noisy. 
The noisiness as well as the overall dynamics resemble the snapshots of the dynamics when modeled by the ABM. 

In Section~\ref{sec:ex_comp}, we will investigate this question further and study a larger sample of simulations as well as different population sizes $N$.

\subsection{Studying the Consistency of the Numerical Solution to the SPDE}\label{sec:num_con}

\begin{figure}[htb!]
\centering
\includegraphics[width=0.8\linewidth]{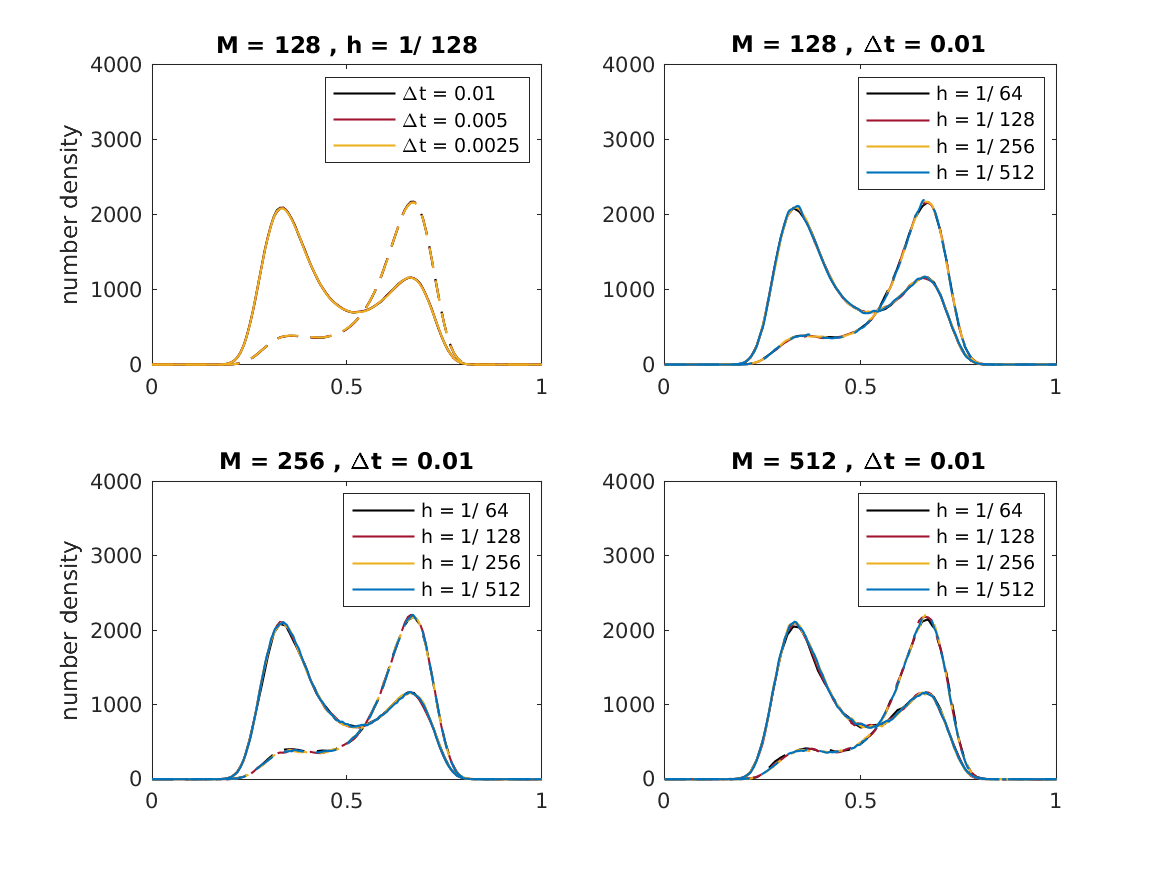}
\caption{For a fixed number of  $N=1000$ agents, of which at   time $t=0$ there are $800$ non-adopters of the innovation (solid line)  and $200$ adopters of the innovation (dashed line), we plot the mean of the discretized solution  over \texttt{sim} $=5000$ realizations  at time $t=1.5$ for different discretization parameters $\Delta t$, $h$, $M$. The remaining parameters are set as in Section \ref{sec:ex_spde}. } \label{fig:cons}
\end{figure}

Before studying the simulated models further, we want to  numerically study the consistency of the numerical solution to the SPDE, i.e., we consider the effect  of varying (a) the time-discretization $\Delta t$, (b) the space-discretization $h$, and (c) the truncation of the noise $M$ for otherwise fixed parameters, with the aim of reaching a fine enough discretization such that the numerical solution to the SPDE is no longer affected by varying the discretization parameters.
 The results  are shown in Figure \ref{fig:cons}. In agreement with these results, we  choose in this chapter for the SPDE discretization $\Delta t= 0.01$, $h=\frac{1}{128}$, and $M=128$.  
 We remark that the choice of the noise truncation threshold $M$ is still an open question, here we choose $M = \frac{1}{h}$ as proposed in \cite{lord2014introduction}.
 
\subsection{Comparison of the Models}\label{sec:ex_comp}
Considering the ABM as the ground-truth model, we compare the innovation spreading dynamics (as introduced in the previous Sections \ref{sec:ex_abm} and \ref{sec:ex_spde}) when  using the ABM and the SPDE model for different population sizes. In particular, we are studying how the discretization of the agent-based and SPDE model compare in computation time and  at which agent numbers $N$ the discretized SPDE model starts to be a good approximation to the simulated ABM. 

We simulate both models for varying population sizes 
$N \in [50,3000]$. For each $N$, we consider an ensemble of \texttt{sim} $=10000$ realizations in order to compute meaningful ensemble averages. The simulation schemes are implemented in Matlab and run on a computer with an Opteron 8384 CPU. 
\subsubsection{The Computational Effort} 
For both model discretizations and varying $N$, we fix the step size in time to $\Delta t=0.01$ and measure the time it takes to simulate one time step. 
Both simulation approaches make use of an Euler-Maruyama time discretization, but space is treated differently. In the ABM approach, we use an Euler-Maruyama discretization for the position dynamics of each agent, whereas in the SPDE approach, we use the scheme for each hat function. Moreover, for the simulation of the ABM, we have to compute  pair-wise distances between agents, which becomes very expensive for increasing agent numbers $N$. Thus we expect the computational effort for the ABM simulation to depend on the number of agents, whereas for the SPDE model it should be independent of $N$ and thus constant for increasing $N$. 

\begin{figure}[htb!]
\centering
\includegraphics[width=0.7\linewidth]{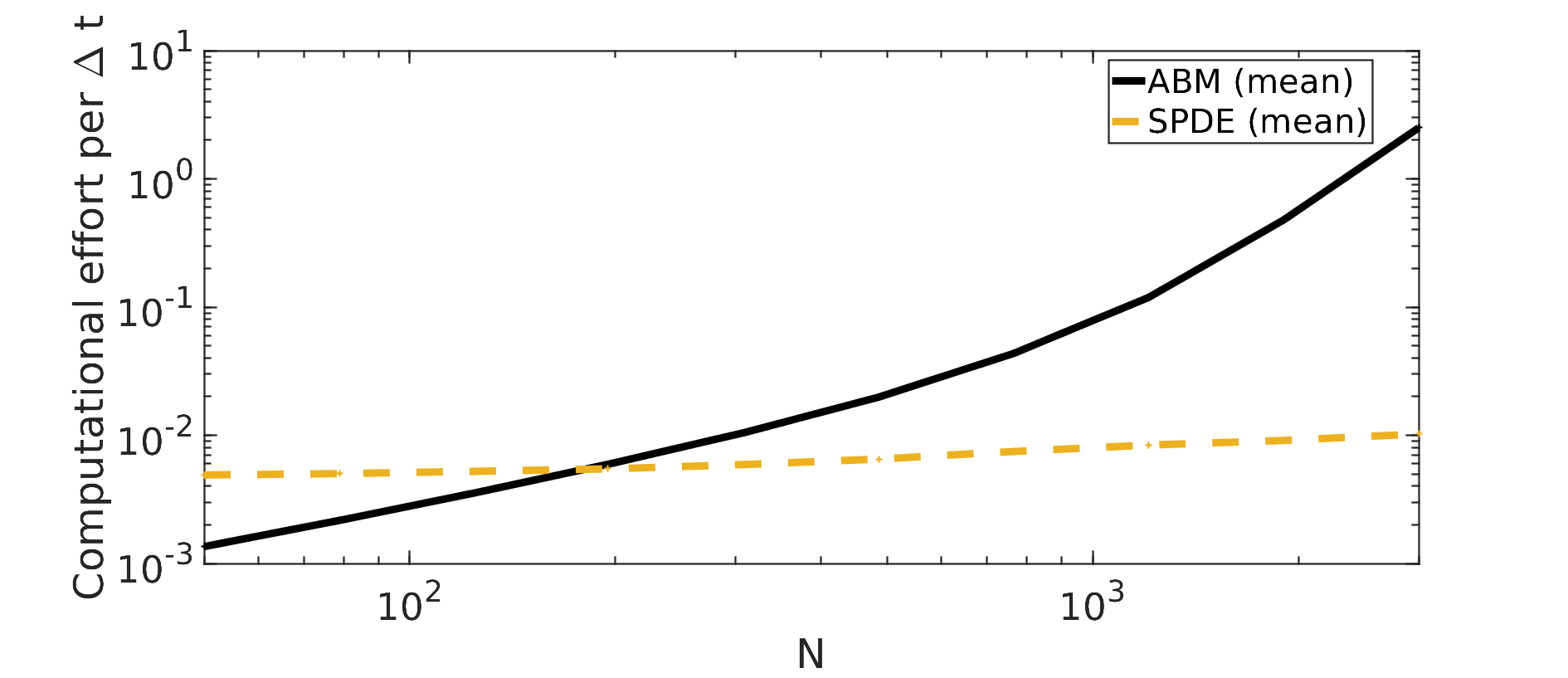}
\caption{The computational effort in s of one time step.}\label{fig:comp1a}
\end{figure}

The results of the computational complexity studies are shown in Figure \ref{fig:comp1a}. As expected, the computational effort for the ABM discretization increases strongly with the number of agents, whereas the effort for the simulation of the density-based model remains cheap. The cost of simulating realizations of the SPDE is independent of $N$ and several magnitudes below the cost of simulating the ABM, except for systems of very few agents, here for approximately $N<10^2$, then the ABM is actually much cheaper to simulate per time step.

\subsubsection{Approximation Quality} 
Next, we compute observables of the simulated dynamics for both models and increasing agent numbers. Based on these measured observables,  we can compare how well they agree and deduce the approximation quality of the SPDE discretization to the ABM discretization.  

In particular we consider:
\begin{itemize}
    \item the time it takes until the  system has reached a certain state for the first time, here we consider the time until $75 \%$  of agents are of type $T_2$.
    We also denote the mean of the first time that  $75 \%$ of agents are adopters in an ABM of $N$ agents by $\tau_N $.
    \item the spatial distribution of type $T_s$ ($s=1,2$) agents  at a fixed model time point, here we choose the time point  $\tau_N$. In order to compare $\tilde{\rho}_s(x,t)$, the agent number density of the discretized SPDE, with the discretized ABM, we first define the empirical agent density $\rho^{\text{ABM}}_s(x,t)$ of the ABM   as the sum of unit masses (e.g.,
a hat function or a narrow Gaussian function that integrate to one) placed at the positions of agents of type $T_s$. 
 Then we can estimate the relative $l_2-$error between the  mean discretized solutions
 \begin{equation}\label{eq:mean_error}
 E_N^{\text{Mean}}:=\frac{\|M(\tilde{\rho}_s(\cdot, \tau_N))-M(\rho^{\text{ABM}}_s(\cdot, \tau_N)) \|_2}{\|  M(\rho^{\text{ABM}}_s(\cdot, \tau_N)) \|_2}
 \end{equation}
for $s=1,2$ and varying $N$ and where $M$ denotes the sample mean, as well as the relative $l_2-$error between
the deviations of the  mean solutions
 \begin{equation}\label{eq:stdev_error}
  E_N^{\text{Std}}:=\frac{\| S( \tilde{\rho}_s(\cdot, \tau_N))-S( \rho^{\text{ABM}}_s(\cdot, \tau_N)) \|_2}{\| S(  \rho^{\text{ABM}}_s(\cdot, \tau_N)) \|_2},
 \end{equation}
 where  $S$ stands for the sample standard deviation.
\end{itemize}
Of course, in computing these observables with the different models we not only make a numerical error due to our numerical discretization, but we also introduce a statistical error when estimating the mean and variances of observables. 

\begin{figure}[htb!]
\centering
\includegraphics[width=0.7\linewidth]{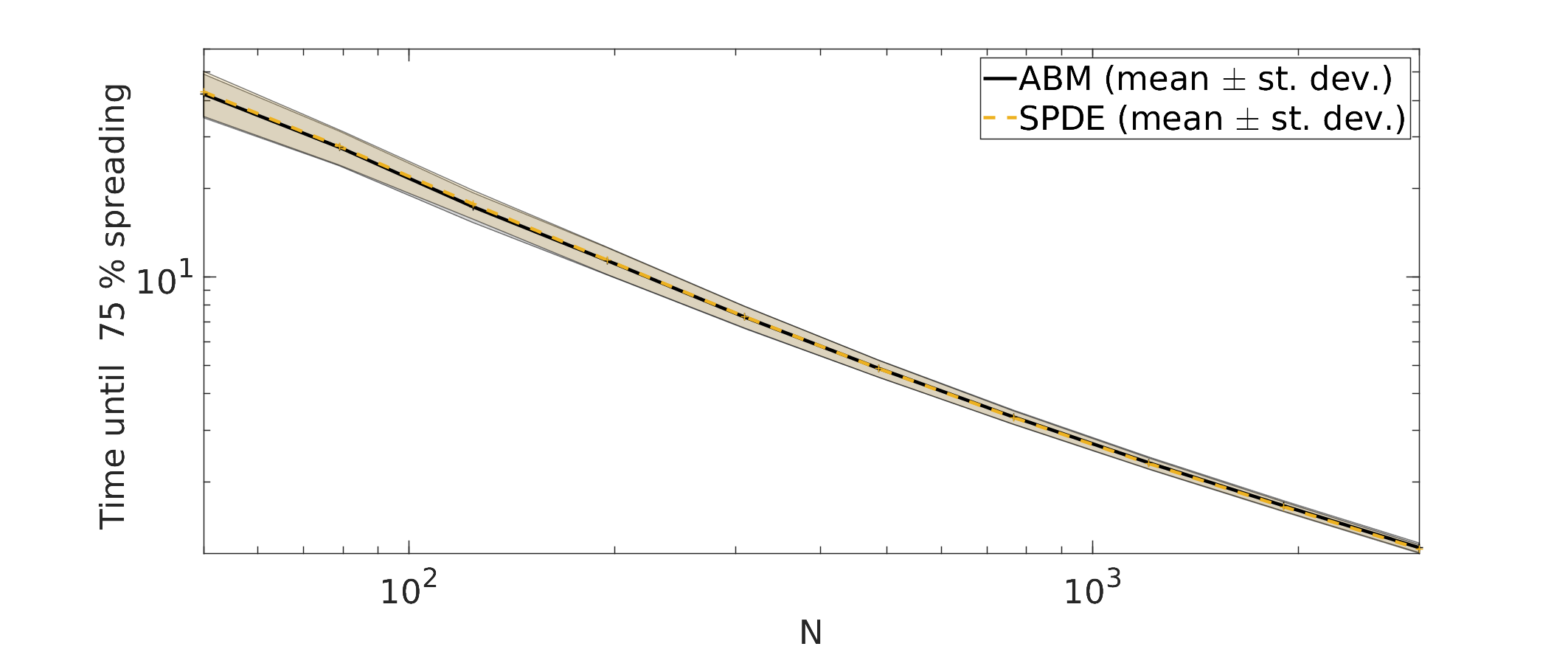}
\caption{ The mean and standard deviation of the first time that $75\%$ of agents are adopters. }\label{fig:comp1b}
\end{figure}

The comparison of the  first time of reaching $75\%$ adopters in the system is shown in  Figure~\ref{fig:comp1b}. The time is very well approximated by the discretized SPDE both in mean and standard deviation when the number of agents is large (in our example around $N>200$), but even for smaller systems the SPDE delivers an acceptable approximation.

\begin{figure}[htb!]
\centering
\includegraphics[width=0.7\linewidth]{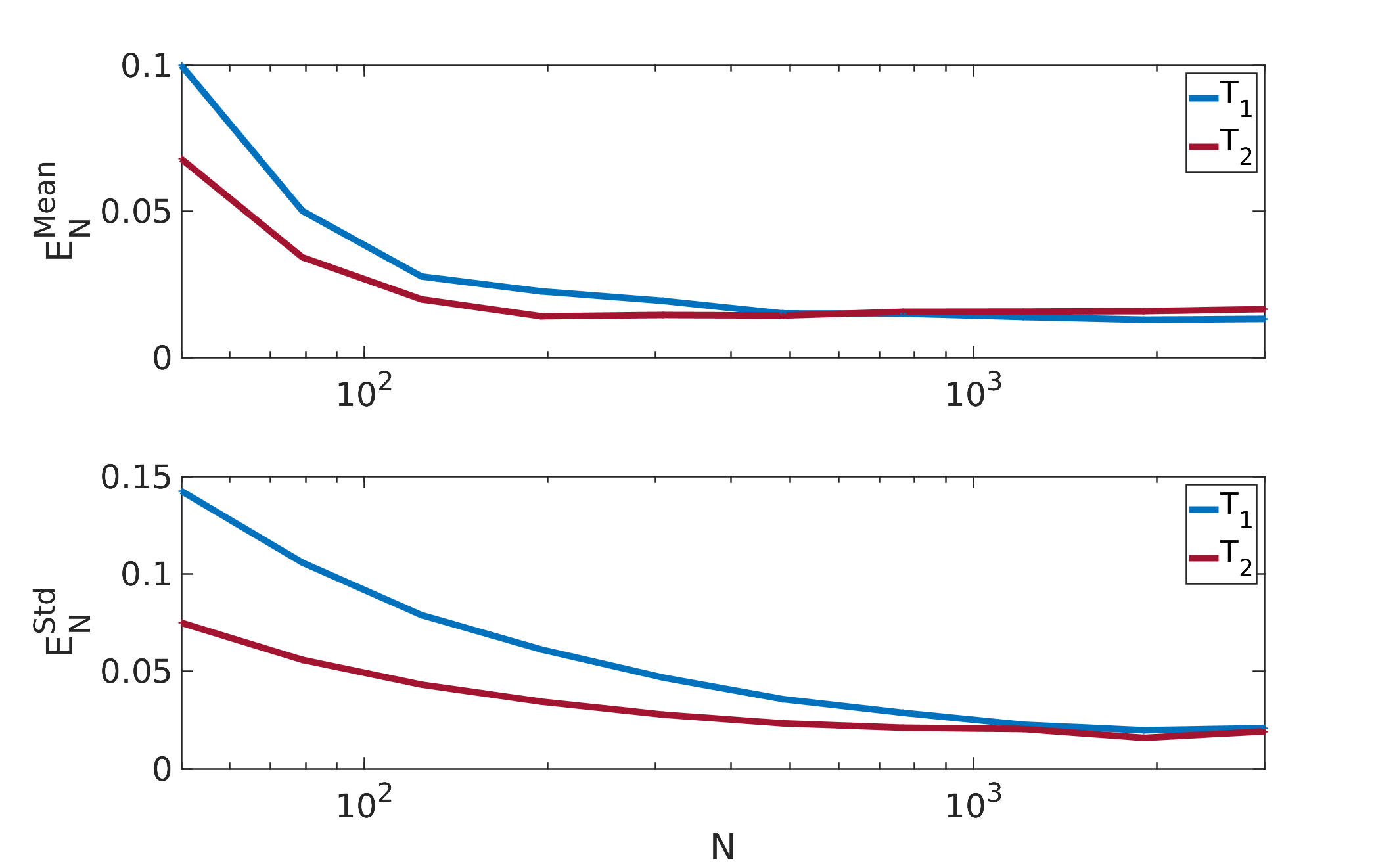}
\caption{ We plot the relative $l_2-$error $E_N^{\text{Mean}}$ between the mean  discretized solutions of the SPDE and the ABM  at some fixed time $\tau_N$ for varying N and also the relative error $E_N^{\text{Std}}$ between the standard deviations of the solutions. }\label{fig:l2}
\end{figure}

\begin{figure}[htb!]
\centering
\includegraphics[width=1\linewidth]{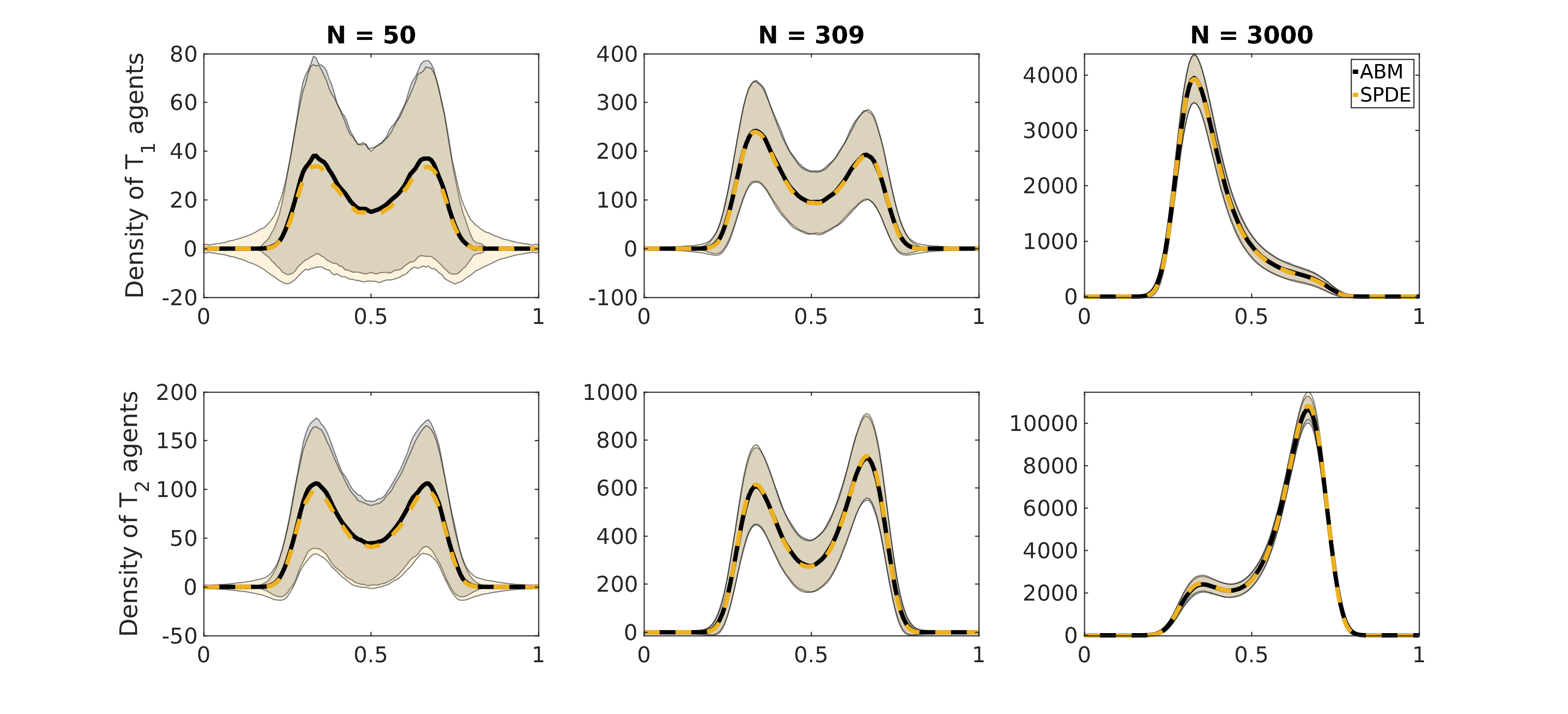}
\caption{The  agent number density (in mean and standard deviation)  for a fixed time point, here taken as $\tau_N$ as given by the black curve in Figure \ref{fig:comp1b}. The outcome  of the ABM is shown in black, the outcome of the SPDE model is depicted in yellow.}\label{fig:comp2b}
\end{figure}

In Figures~\ref{fig:l2} and~\ref{fig:comp2b} we study the spatial distribution of agents at the fixed time point $\tau_N$ for increasing population sizes $N$. From Figure~\ref{fig:comp2b} we can note that for small $N=50$, the standard deviation of the  number density near the domain boundaries $x=0$ and $x=1$ is not that well approximated due to the small number of agents.  
Figure~\ref{fig:l2} indicates that the error  between the solution (in mean and standard deviation) of the discretized ABM and SPDE decreases with increasing $N$.
Still we do not observe a convergence of the errors to zero in the large population limit, which can be explained as the consequence of the following two error sources:
First, for different $N$, we are comparing qualitatively different agent densities. For small $N$, the agent number densities are rather spatially equilibrated at time point $\tau_N$ (i.e., they resemble the Boltzmann distribution), while for large $N$, the agent densities are spatially not equilibrated at time $\tau_N$ (see Figure~\ref{fig:comp2b}). This is due to the spatial diffusion of agents being independent of the number of agents, while the innovation diffusion is much faster when there are more agents in the system. 
Second, we are not comparing the exact solutions of ABM and SPDE in the large population limit, but numerical solutions, thus both solutions have a numerical error that does not vanish in the large population limit.

In general we can remark that for increasing numbers of agents, the stochasticity in both simulation methods decreases as can be seen from the decrease in standard deviation. The more agents in the ABM, the less the individual behaviour sticks out (due to the law of large numbers).  For the SPDE model this was also expected (see the comment in Section~\ref{sec:complete}), since the noise scaling of the SPDE is such that the  noise dominates for small population sizes and has less influence for large population sizes.

\section{Conclusion and Future Outlook}
Agent-based models are  the most natural approach to construct models for real-world systems containing discrete interacting entities such as humans, since one  simply has to describe  the actions of each individual. 
Based on~\cite{conrad2018mathematical, ourepjpaper}, we have introduced a general agent-based model that is formulated in terms of coupled diffusion and Markov jump processes for each agent.

Some further research regarding the construction of the ABM could be directed at the following questions:
     \begin{itemize}
         \item How reasonable is it to model human or agent mobility by drifted Brownian motion or should we  instead describe the diffusion of agents on an infrastructure network or by a different stochastic process~\cite{barbosa2018human}?
         \item Can we extend the ABM to include feedback loops, such that e.g., the type changes also influence the position dynamics?
     \end{itemize}
Simulating agent-based models for real-world dynamics quickly becomes costly due to an explosion in the computational complexity for increasing agent numbers and the need for repeated simulations due to its stochastic description. For instance, when modeling the spreading of innovations in ancient times such as in~\cite{conrad2018mathematical,ourepjpaper}, many Monte-Carlo simulation are required to capture the full spectrum of the diverse dynamics. But this becomes computationally very expensive. A thorough sensitivity analysis of the parameters demands many simulations for each parameter set and is as such not tractable.

Based on an extension of the Dean-Kawasaki model~\cite{kim2017stochastic,dean1996langevin,kawasaki1973new,bhattacharjee2015fluctuating,helfmann2018stochastic}, we therefore considered a meso-scale approximation to the ABM for systems of many agents. The reduced model is given by a system of coupled stochastic PDEs propagating agent densities for the different agent types.

For both models, the ABM and the reduced density-based model, we constructed and explained simulation schemes. The Finite Element discretization of the system of SPDEs serves as an interpretation and regularization of the  ill-defined SPDE for which solutions not necessarily exist~\cite{fehrmanwell,cornalba2018regularised,lehmann2018dean}.

There are however several questions that still remain unanswered and should be investigated further:
     \begin{itemize}
         \item The mapping from the ABM to the SPDE for more complex interaction rules and dynamics has to be worked out in detail.
         \item The notion of the solution to the SPDE and its existence and uniqueness  is still unclear due to several mathematical difficulties ~\cite{cornalba2018regularised,lehmann2018dean,fehrmanwell}.
         \item Can we find quantitative statements for the approximation (especially in different limits, e.g., small or large populations,  diffusion speeds and interaction rates) and discretization quality?
         \item Numerical experiments for higher spatial dimensions and more complex domains as well as studies on the applicability to real-world systems containing humans need to be investigated further.
     \end{itemize}
Finally, we compared the simulation effort and studied the approximation quality of the reduced density-based model to the ABM computationally on a toy example of innovation spreading. From these computational experiments we can conclude the following. For systems of few agents, the SPDE approach is more costly than simulating the ABM. We consider the ABM as the ground-truth model and thus as the most accurate description of the agent system. The dynamics have to be described in terms of individual agents, since there are only very few. 

But for systems of many agents, we can instead use the approximation by the SPDE to study the agent system dynamics much more efficiently. In the case of our toy example, this approximation is very accurate in mean and standard deviation already for systems of a few hundred agents.  Thus, the presented reduced model is a very promising tool for modeling and especially simulating real-world systems with large agent populations.
The SPDE model offers the possibility to carry out a more thorough model analysis such as parameter studies and inference, or model control, whilst staying numerically tractable.

\paragraph{Acknowledgments}
The authors thank Changho Kim  for an insightful discussion on the Dean-Kawasaki equation with added reactions, and the anonymous reviewers for detailed and helpful feedback.
The authors received funding by the Deutsche Forschungsgemeinschaft (DFG, German Research Foundation) under Germany's Excellence Strategy – The Berlin Mathematics
Research Center MATH+ (EXC-2046/1, project ID: 390685689) and through grant CRC 1114.

\bibliographystyle{alpha}.
\small{\bibliography{references}}
%\bibliographystyle{unsrt}
%\bibliography{references}
\end{document}